\documentclass{elsart}

\usepackage[numbers, sort&compress]{natbib}
\usepackage[utf8]{inputenc}
\usepackage{amsmath,amssymb,amstext}

\usepackage{pifont}
\usepackage{algorithm}
\usepackage[noend]{algorithmic}
\usepackage{supertabular}
\usepackage{latexsym,epic,eepic,epsfig,graphicx,euscript}
\usepackage{psfrag}

\usepackage{subfigure}
\usepackage{xspace}
\usepackage{verbatim} 

\renewcommand{\algorithmicensure}{\textbf{Output:}}
\renewcommand{\algorithmicrequire}{\textbf{Input:}}

\providecommand{\figref}[1]{Figure~\ref{#1}}
\providecommand{\Figref}[1]{Figure~\ref{#1}}
\providecommand{\secref}[1]{section~\ref{#1}}

\providecommand{\probref}[1]{Problem~\ref{#1}}

\providecommand{\tabref}[1]{Table~\ref{#1}}
\providecommand{\Tabref}[1]{Table~\ref{#1}}
\providecommand{\eqnref}[1]{equation~\eqref{#1}}

\providecommand{\thmref}[1]{Theorem~\ref{#1}}

\providecommand{\remref}[1]{Remark~\ref{#1}}

\providecommand{\propref}[1]{Proposition~\ref{#1}}

\providecommand{\algref}[1]{Algorithm~\ref{#1}}
\providecommand{\Algref}[1]{Algorithm~\ref{#1}}
\providecommand{\defref}[1]{Definition~\ref{#1}}

\providecommand{\lemref}[1]{Lemma~\ref{#1}}

\providecommand{\corref}[1]{Corollary~\ref{#1}}

\newcommand{\N}{\mathbb{N}}
\newcommand{\Z}{\mathbb{Z}}

\newcommand{\K}{\mathbb{K}}

\newcommand{\Ztwo}{\Z_2}
\newcommand{\ZTwo}{\Ztwo}
\newcommand{\cR}{C}                 
\newcommand{\wR}{R}                 
\newcommand{\primes}{\mathbf{p}}    

\newcommand{\calG}{\mathcal{G}}     %
\newcommand{\vars}{\mathbf{x}}      
\newcommand{\monom}{\mathbf{x}}     

\DeclareMathOperator{\ecart}{ecart}       %
\DeclareMathOperator{\image}{Im}          %
\DeclareMathOperator{\sPoly}{spoly}       %
\DeclareMathOperator{\lt}{LT}             %
\DeclareMathOperator{\lc}{LC}             %
\DeclareMathOperator{\lexp}{LE}           %
\DeclareMathOperator{\lm}{LM}             %
\DeclareMathOperator{\li}{L}              %
\DeclareMathOperator{\lmi}{LM}            %
\DeclareMathOperator{\nf}{NF}             %
\DeclareMathOperator{\lcm}{lcm}           %
\DeclareMathOperator{\gcod}{gcd}          %
\DeclareMathOperator{\NEid}{NE}           
\DeclareMathOperator{\Nid}{N}             
\DeclareMathOperator{\syz}{Syz}           
\DeclareMathOperator{\Crit}{crit}         
\DeclareMathOperator{\smaranch}{SM}       
\DeclareMathOperator{\ZeroDiv}{NT}        

\newcommand{\skalar}[1]{\left\langle #1\right\rangle}  
\newcommand{\card}[1]{{\left\lvert #1\right\rvert}}    
\newcommand{\Image}[1]{\image\left( #1\right)}         
\newcommand{\Syz}[1]{\syz\left( #1\right)}             
\newcommand{\wo}{\backslash}                           
\newcommand{\units}[1]{{#1}^*}           
\newcommand{\zeros}[1]{\Nid\left( #1\right)}           
\newcommand{\nonunits}[1]{\NEid\left( #1\right)}       
\newcommand{\LT}[1]{\lt\left( #1\right)}               
\newcommand{\LM}[1]{\lm\left( #1\right)}               
\newcommand{\LC}[1]{\lc\left( #1\right)}               
\newcommand{\LE}[1]{\lexp\left( #1\right)}             
\newcommand{\LI}[1]{\li\left( #1\right)}               
\newcommand{\LMI}[1]{\lmi\left( #1\right)}             
\newcommand{\kgv}[1]{\lcm\left( #1\right)}             
\newcommand{\ggt}[1]{\gcod\left( #1\right)}            
\newcommand{\spolyr}[1]{\sPoly_r\left( #1\right)}      
\newcommand{\zeroDiv}[1]{\ZeroDiv\left( #1\right)}     

\newcommand{\NF}[2]{\nf\left(\left. #1\,\right|\,#2\right)} 

\newtheorem{notation}[thm]{Notation}
\newtheorem{example}[thm]{Example}
\newtheorem{remark}[thm]{Remark}
\newtheorem{problem}[thm]{Problem}

\newtheorem{observation}[thm]{Observation}

\newtheorem{theorem}[thm]{Theorem}
\renewcommand{\algorithmiccomment}[1]{/* #1 */}
\newcommand{\set}[1]{\{#1\}}

\newcommand{\Singular}{{\sc Singular}\xspace}
\newcommand{\Sage}{SAGE\xspace}
\newcommand{\Bool}{\ensuremath{\mathbb{B}}\xspace}

\newtheorem{lemma}[thm]{Lemma}
\newtheorem{proposition}[thm]{Proposition}
\newtheorem{corollary}[thm]{Corollary}

\newtheorem{definition}[thm]{Definition}

\DeclareMathOperator{\tail}{tail}
\DeclareMathOperator{\spoly}{spoly}
\DeclareMathOperator{\varsof}{vars}
\DeclareMathOperator{\V}{V}
\DeclareMathOperator{\I}{I}

\DeclareMathOperator{\support}{supp}

\DeclareMathOperator{\True}{True}
\DeclareMathOperator{\BI}{BI}
\DeclareMathOperator{\BGB}{BGB}

\DeclareMathOperator{\symmgbGFTwoOp}{symmgbGF2}
\newcommand{\symmgbGFTwo}{\ensuremath{\symmgbGFTwoOp}\xspace}

\DeclareMathOperator{\thenBranch}{\text{\texttt{then}}}
\DeclareMathOperator{\elseBranch}{\text{\texttt{else}}}
\DeclareMathOperator{\blockdeg}{blockdeg}
\DeclareMathOperator{\ite}{ite}
\DeclareMathOperator{\successor}{succ}

\newcommand{\ie}{i.\,e.\xspace}
\newcommand{\eg}{e.\,g.\xspace}
\newcommand{\sth}{s.\,th.\xspace}
\newcommand{\mcomma}{\,\mbox{,}}
\newcommand{\mfstop}{\,\mbox{.}}

\newcommand{\defemph}{\emph}

\newcommand{\Magma}{{\sc Magma}}

\newcommand{\PolyBoRi}{{\sc PolyBoRi}\xspace}

\newcommand{\explZtwoxoneton}{{\Ztwo[x_1, \dots, x_n]}}
\newcommand{\Ztwoxoneton}{{\ensuremath{\Ztwo[\mathbf{x}]}}\xspace}

\newcommand{\FP}{\ensuremath{\mathsf{FP}}\xspace}
\newcommand{\fieldequations}{\FP}
\newcommand{\explfieldequations}{{x_1^2+x_1,\dots,x_n^2+x_n}}
\newcommand{\fieldequationset}{\FP}
\newcommand{\ideal}[1]{\ensuremath{\langle #1\rangle}\xspace}

\newcommand{\fieldideal}{\ideal{\FP}}

\newcommand{\pwithfe}[1]{\ideal{#1, \fieldequations}}

\newcommand{\wrt}{w.\,r.\,t.\xspace}

\newcommand{\behaviour}{behaviour\xspace}
\newcommand{\optimis}{optimis}
\newcommand{\optimised}{\optimis{}ed\xspace}

\newcommand{\optimisation}{\optimis{}ation\xspace}
\newcommand{\optimisations}{\optimisation{}s\xspace}
\newcommand{\specialised}{specialised\xspace}
\newcommand{\modell}{modell}
\newcommand{\modelling}{\modell{}ing\xspace}

\newcommand{\summarised}{summarised\xspace}
\newcommand{\minimise}{minimise\xspace}

\newcommand{\utilises}{utilises\xspace}
\newcommand{\utilised}{utilised\xspace}

\newcommand{\utilising}{utilising\xspace}
\newcommand{\homogenisation}{homogenisation\xspace}
\newcommand{\generalisation}{generalisation\xspace}

\newcommand{\emphasise}{emphasise\xspace}

\providecommand{\color}[1]{}

\newcommand{\powerset}[1]{\ensuremath{\mathcal{P}(#1)}}

\newcommand{\Groebner}{Gr\"{o}bner\xspace}
\newcommand{\naive}{na{\"i}ve\xspace}

\renewenvironment{proof}{\begin{pf}}{\end{pf}}

\providecommand{\boolemult}{\ensuremath{{\star\hspace{.075em}}}\xspace}

\begin{document}

\begin{frontmatter}

\title{New developments in the theory of \Groebner bases and
applications to formal verification}

\author{Michael Brickenstein\thanksref{thxmfo}},
\ead{brickenstein@mfo.de}
\author{Alexander Dreyer\thanksref{thxitwm}},
\ead{alexander.dreyer@itwm.fraunhofer.de}
\author{Gert-Martin Greuel\thanksref{thxunikl}},
\ead{greuel@mathematik.uni-kl.de}
\author{Markus Wedler\thanksref{thxunikl}},
\ead{wedler@eit.uni-kl.de}
\author{Oliver Wienand\thanksref{thxunikl}}
\ead{wienand@mathematik.uni-kl.de}

\thanks[thxmfo]{Mathematisches Forschungsinstitut Oberwolfach,
Schwarzwaldstr.\ 9-11, 77709 Oberwolfach-Walke, Germany}

\thanks[thxitwm]{Fraunhofer Institute for Industrial Mathematics~(ITWM)\\
Fraunhofer-Platz 1, 67663 Kaiserslautern, Germany}

\thanks[thxunikl]{University of Kaiserslautern,
Erwin-Schr\"{o}dinger-Stra{\ss}e, 67653 Kaiserslautern, Germany}

\begin{abstract}
We present foundational work on standard bases over rings and on Boolean
\Groebner bases in the framework of Boolean functions. The research was
motivated by our collaboration with electrical engineers and computer scientists
on problems arising
from formal verification of digital circuits. In fact, algebraic \modelling of
formal verification problems is developed on the word-level as well as on the
bit-level. The word-level model leads to \Groebner basis in the polynomial
ring over~$\Z/2^n$ while the bit-level model leads to Boolean \Groebner
bases. In addition to the theoretical foundations of both approaches, the
algorithms have been implemented. Using these implementations we show that
special data structures and the exploitation of symmetries make \Groebner bases
competitive to state-of-the-art tools from formal verification but having
the advantage of being systematic and more flexible.
\end{abstract}

\begin{keyword}
\Groebner basis, formal verification, property checking, Boolean polynomials,
satisfiability
\end{keyword}

\end{frontmatter}

\section*{Introduction}
It has become common knowledge in many parts of mathematics and in some
neighbouring fields that \Groebner bases are a universal tool for any kind of
problem which can be modelled by polynomial equations. However, quite often the
models involve too many unknowns and equations making it unfeasible to carry
out the corresponding \Groebner basis computation.

This is, for example, the case for most real-world problems from discrete
\optimisation or from formal verification of digital systems, two areas of
eminent practical importance. Because of their importance the community
working in these fields is much bigger than the \Groebner basis community
and, moreover, there exist highly \specialised commercial tools making it
unrealistic to believe that \Groebner bases can be of comparable practical
efficiency in these areas.

One of the purposes of this paper is to show that, in many cases \Groebner
bases can be used to find solutions for formal verification problems.
In this way, this forms a good complement to existing
techniques, like simulators and SAT-solver, which are suited for identification
of counter examples~(falsification).

A significant advantage is, that \Groebner bases
provide a mathematically prov\-en systematic and very flexible tool while many
engineering solutions inside commercial verification tools rely on ad hoc
heuristics for special cases. However, the success of \Groebner basis methods,
reported in this paper, could not be achieved with existing generic \Groebner
basis algorithms and implementations. On the contrary, it relies on the theory
of \Groebner bases in Boolean rings and improvements of algorithms for this
case, both being developed by the authors and described here for the first
time.

The Boolean \Groebner basis formulation of a verification problem comes from a
modelling on the bit-level. We describe here also another approach based on a
modelling on the word-level, leading to \Groebner basis computations in the
polynomial ring over the ring $\Z_{2^n}$ of integers modulo $2^n$ where $n$ is
the word length, that is, the number of bits used by each signal. This
approach has the advantage that it leads to a more compact formulation with
less variables and equations. On the other hand, it has the disadvantage that
$\Z_{2^n}$ is not a field for $n>1$, but a ring with zero divisors. Moreover, we show
that an arbitrary verification problem cannot, in general, be modelled by a
system of polynomial equations over the ring $\Z_{2^n}$ and, furthermore, we can
in general only prove non-satisfiability but not satisfiability. Nevertheless,
a combination of the word-level with the bit-level model could overcome these
difficulties by preserving some of the advantages of the word-level
approach. However, this is not yet fully explored and hence not presented in
this paper.

The paper is organized as follows. In section 1 we describe the formal
verification of digital circuits and its algebraic modelling via word-level
and bit-level encoding. We do also discuss the advantages and disadvantages
of both approaches.

The second section presents foundational results about standard bases in
polynomial rings over arbitrary rings, allowing monomial orderings which are
not well orderings. New normal form algorithms and criteria for $s$-polynomials
 are presented in the case of weakly factorial principal ideal
rings. This includes the case $\Z_m$ which is of interest in the application to
formal verification.

In section 3 the theory of Boolean \Groebner bases is developed in the
framework of Boolean functions. Mathematically the ring of Boolean functions~%
$\Z^n_2\to \Z_2$ is isomorphic to $\explZtwoxoneton/\fieldideal$
where $\FP$ is the set of field polynomials~$x_i^2+x_i$, for~$i=1,\ldots, n$. Boolean
\Groebner bases are \Groebner bases of ideals in~$\Ztwoxoneton$
  containing~$\FP$, modulo the ideal~$\fieldideal$. The usual data
  structure for polynomials in~$\Ztwoxoneton$ is, however, not
  adequate.

We propose to encode Boolean polynomials as zero-suppressed binary decision
diagrams (ZDDs) and describe the necessary algorithms for polynomial arithmetic
which takes advantage of the ZDD data structures. Besides the polynomial
arithmetic the whole environment for \Groebner basis computations has to be
developed. In particular, we describe
efficient comparison algorithms for the
most important monomial orderings. A central observation, which is responsible
for the success of our approach (besides the efficient handling of the new
data structures), is the appearance of \emph{symmetries} in systems of Boolean
polynomials coming from formal verification. The notion of a symmetric monomial
ordering is introduced and an algorithm making use of the symmetry is
presented.

The presented algorithms have all been implemented, either in \Singular
or in the \PolyBoRi-framework.

In the last chapter we present some implementation details and explicit
timings, comparing the new algorithms with state-of-the-art implementations
of either \Groebner basis algorithms or SAT-solvers. Moreover, we discuss
open problems, in particular for polynomial systems over $\Z_{2^n}$.

\section*{Acknowledgements}

The present research is supported by the Deutsche Forschungsgemeinschaft within
the interdisciplinary project ``Entwicklung, Implementierung und Anwendung
mathematisch-algebraischer Algorithmen bei der formalen Verifikation
digitaler Systeme mit Arithmetikbl\"ocken'' together with the research group
of Prof.\ W.\ Kunz from the department ``Electrical and Computer Engineering'' at
the University of Kaiserslautern.

Moreover, the work was also supported by the Cluster of Excellence in
Rhine\-land-Palatinate within the DASMOD and VES projects. We like to thank
all institutions for their support.

This paper is an enlarged version of a talk by the third author given at the
RIMS International Conference on ``Theoretical Effectivity and Practical
Effectivity of \Groebner Bases'' in Kyoto, January 2007. We like to thank
T.\ Hibi for organizing this conference and for his hospitality.

\section{Algebraic models for formal verification}
\subsection{Formal verification}

The presented research was spurred by a joint project on formal
verification with the electrical engineering department at the University of
Kaiserslautern.
An important goal pursued in modern circuit design flows is to avoid
the introduction of bugs into the circuit design in every stage of the
process.
We do not go into detail here, but just mention, that
formal verification of hard- and software is a huge field of research with an overwhelming
amount of literature. We refer to \citep{McMillan, HachtelBook, kunzetal01} for more details and references.

Property checking is a technique for functional
verification of the initial register transfer level (RTL) description
of a circuit design. The initial specification of the design that is
often given as a more or less informal human readable document is
formalized by a set of properties. A systematic methodology ensures
that the complete intended behavior of the circuit is covered by the
resulting property suite. However, each property describes the
required circuit behavior in a well defined scenario. This allows for
an early evaluation for parts of the design as soon as they are
completed. 


Classical methods for design validation include the simulation of the system
with respect to suitable input stimuli, as well as, tests based on emulations,
which may use simplified prototypes.  The latter may be constructed using field
programmable gate arrays~(FPGAs).
Due to a large number of possible settings, these approaches can never cover the
overall \behaviour of a proposed implementation.  In the worst case, a defective
system is manufactured and delivered, which might result in a major product
recall and liability issues.
Therefore simulation methods are more and more replaced by formal methods which
are based on \emph{exact} logical and mathematical
algorithms for automated proving of circuit properties.

\subsection{Design flow}

\begin{figure}
\begin{center}
\includegraphics[width=.7\textwidth]{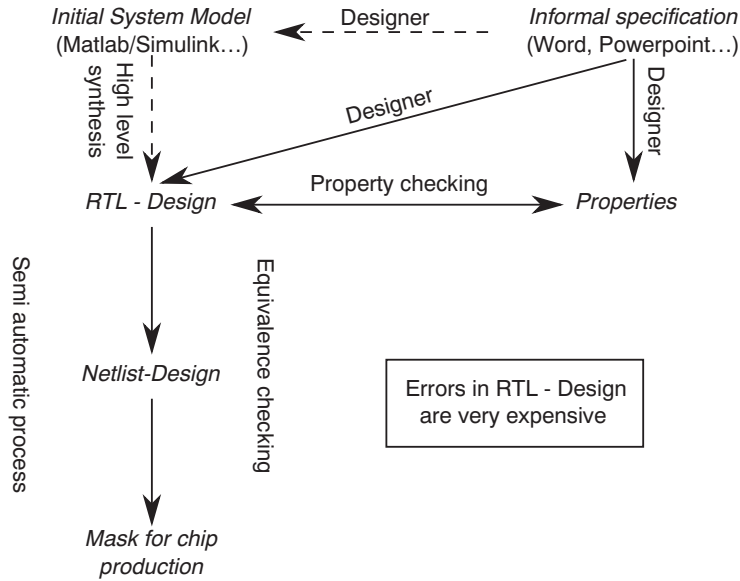}
\end{center}
\caption{Digital system design flow\label{fig:production}}
\end{figure}

The circuit design starts with an informal specification of a
microchip~(\figref{fig:production})
by some tender documents which are usually given in a human readable text or
presentation format. In a first step the specification may be translated in a
highlevel \modell{}ing language.
%
%
One possibility is to use
high level synthesis for generating a \defemph{register transfer level}~(RTL)
design which describes the flow of signals between registers in terms of
a hardware description language~\cite{vhdlverilog}.
But this is rarely used in practise as it does constrain the freedom of the
design.
Instead, designers manually create the RTL design in a hardware desription language .
Concurrently, intended behavior specified by the informal specification is
formalized by formal properties.
Automatic tools are used to ensure that the
RTL design  fulfills these conditions.

%
After passing property checking a netlist is generated semi-automatically
from the RTL. The latter is used to derive the actual layout of  the chip mask.
The validation that different circuit descriptions arising from  the last two
steps  emit the same \behaviour, is called
\defemph{equivalence checking}. Since this can be handled accurately, setting of
the RTL design is the most crucial part.
Errors at this level may become very expensive, as they may lead to unusable
chip masks or even defective prototypes.
The present paper is concerned with this critical level.

The ability of checking the validity of a proposed design restricts the design
itself: a newly introduced design approach may not be used for an implementation
as long as its verification cannot be ensured.
In particular, this applies to digital systems consisting of combined logic and
arithmetic blocks, which may not be treated with \specialised approaches.
Here, dedicated methods from
computer algebra
may lead to more generic procedures, which help to fill the design gap.

\subsection{Problem formulation and encoding in algebra}

The verification problem is defined by a set of axioms~$M$ representing the
circuit \wrt given decision variables. In addition, a set of statements~$P$
represents the property to be checked.
For instance, if~$M$ models a multiplication unit,
a suitable~$P$ would be the condition that after a complete
cycle the output of~$M$ is the product of its inputs.

The question, whether the circuit represented by~$M$
fulfills~$P$ can be reformulated in the following way:
First of all, we may assume, that~$M$ is consistent, \ie there are no
contradictions inherent in the axioms, since the axioms describe a circuit.
Then the new set of axioms~$M\wedge \neg P$ is contradictable if and only if~$M$
implies~$P$. Hence the desired property $P$ will be proven by showing, that~$M\wedge \neg P$
has no valid instance, \ie one fulfilling the axioms and not the property.

In the following we encode this logical system into a system of algebraic equations in two ways,
on word-level and on bit-level. The word-level model will lead to consider \Groebner bases over
the ring~$\Z_{2^n}$ while the bit-level will lead to \Groebner basis over Boolean rings.
Here and in the following~$\Z_m$ denotes the finte ring~$\Z/m\Z$ for~$m\in\Z\wo\{0\}$.

\subsubsection{Word-level encoding}
\newcommand{\Ztwon}{\ensuremath{\Z_{2^{n}}}}
\label{sec:wordlevel}

\begin{figure}
\begin{center}
\subfigure[RTL diagram]{%
\includegraphics[height=.1\textheight]{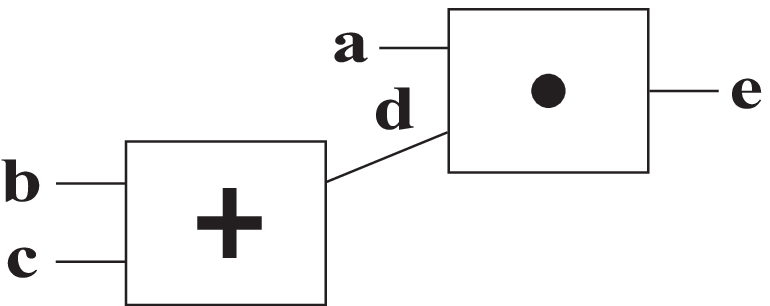}\label{fig:RTLdiagram}
}%
\subfigure[Property]{%
\hspace{1.5em}%
\includegraphics[height=.1\textheight]{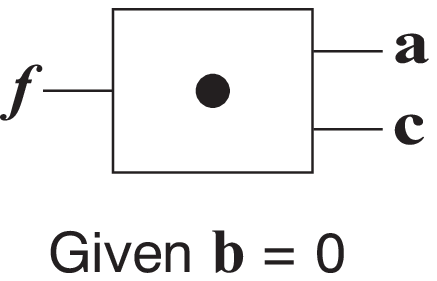}\label{fig:RTLproperty}
}%
\end{center}
\caption{RTL design and property
\label{fig:RTL_design}}
\end{figure}

We illustrate, how the problem of formal verification can be encoded in a system of algebraic
equations using polynomials over the ring~$\Z_{2^n}$.
Let~$n$ be the word length of the circuit, \ie the number
of bits used by each signal (in typical applications we have~$n\in\{16,32,64\}$).
Then the RTL description displayed in \figref{fig:RTLdiagram}  is
equivalent to the following set of algebraic equations
\begin{equation}
\label{eqn:RTLdiagram}
M=\{b + c = d, \; a\cdot d = e\}
\end{equation}
where $b + c - d, a\cdot d - e$ are polynomials in $\Ztwon[a,b,c,d,e,f]$. Of course,
the two equations in~$M$ are equivalent to $a\cdot(b+c)=e$, but in general the latter
input-output form is infeasible due to its complexity. Also, there can be more
than one output per block and only some of these outputs may be used further.

For example, \figref{fig:RTLproperty} presents the property
\begin{equation}
\label{eqn:RTLproperty}
P=\{b=0,\; a \cdot c = f\}.
\end{equation}
In this case, the statement that~$M$ implies~$P$ is equivalent to the assertion
that $M \cup P \cup \{f  \neq e\}$ has no solution.
Since the set $\{f  \neq e\}$ is not a closed algebraic set, we
replace~$f \neq e$ by~$s\cdot(f-e)=2^{n-1}$, where $s$ is a new variable. Indeed, it is easy to see that a value~$s\in\Ztwon$
fulfills this equation if and only if~$f \neq e$ (since the ring~$\Ztwon$ has zero-divisors,~$f\neq e$ cannot
be encoded by~$s(f-e)=1$). Let~$I$ be the ideal
$\skalar{\set{b + c - d,  a\cdot d - e,\;\; b, a \cdot c - f,\;\;s\cdot(f-e)-2^{n-1}}}$
in~$\Ztwon[a,b,c,d,e,f,s]$. Then the question reduces to the question whether
$$
\V(I):=\set{(a,b,c,d,e,f,s)\in\Ztwon^7 \,|\,
    p(a,b,c,d,e,f,s)= 0,\; \text{for all}\; p\in I}$$
is empty. There are no solutions for the ideal $I$ (\ie $\V(I) = \emptyset$) if and only if $M\wedge\neg P$ is
contradictable, that is,~$P$ is satisfied by~$M$.

One way of tackling this problem is to compute a \Groebner basis of~$I$ in
the ring~$R/I_0$, where~$I_0$ denotes the ideal of vanishing polynomials in~$R$, \ie
polynomials evaluating to zero at any point of $\Ztwon^7$. Due to the zerodivisors in this
ring the ideal $I_0$ has more structure than in the finite field case and even its \Groebner basis can become huge~(cf.~\cite{vangb-oli}).

\subsubsection{Bit-level encoding}
\label{sec:bitlevel}                      
An alternative approach is to encode the problem at the bit-level, that is, as polynomials
over~$\Ztwo$.
This approach is based on the fact that every value of~$x$ in~$\Z_{2^{n}}$ can be encoded uniquely to the base 2, \ie in its bits:
\begin{equation}\label{eqn:word2nbits}
x = x_0 + x_1 2 + \dots + x_{n-1} 2^{n-1},\,x_i\in\{0,1\}\mfstop
\end{equation}
In the example above we can express each variable $a,b,c,d,e,f$ analogously to \eqnref{eqn:word2nbits}
with new variables~$a_i, b_i, c_i, d_i, e_i, f_i\in\{0,1\}, i=0,\dots,n-1$. Then
\eqnref{eqn:RTLdiagram} and~\eqnref{eqn:RTLproperty} must be rewritten, which
yields~$n$ equations for each of them.
Gathering all corresponding polynomials and
adding the polynomial~$\prod \left(1 - f_i + e_i\right)$, which is equivalent to~$f\neq e$,
we obtain an ideal~$I$ over~$R:=\ZTwo[a_0, \ldots, f_{n-1}]$ in~$6\,n$ variables.

For instance, the bits~$p_0,\dots,p_{n-1}\in\{0,1\}$
of the product~$p=a\cdot b$ are given by
equations~$p_j=a_{j}\cdot b_{0} + \sum_{i=0}^{j-1} (a_{i}\cdot b_{j-i} +  t_{i,
j-i})$ over~$\Ztwo$, where the~$t_{k,l}$ mark rather complicated bit-level
expressions in the~$s_{k,l}\in\{0,1\}$,
which fulfill~$ p_{k} + s_{k,1} 2 + \dots + s_{k,n-1} 2^{n-1} =
       a_{k}\cdot b_{0} + \sum_{i=0}^{k-1} (a_{i}\cdot b_{k-i} +  s_{i,k-i}) $
in~$\Z_{2^{n}}$.
For example, for~$n=4$, we get
\begin{eqnarray*}
p_3 &=& a_3\,b_0 + a_2\,b_1+ a_1\,b_2 + a_0\,b_3
 + a_2\,a_1\,a_0\,b_1\,b_0 +\\&& \quad a_2\,a_1\,b_1\,b_0 +
a_2\,a_0\,b_2\,b_0  + a_1\,a_0\,b_2\,b_1\,b_0 + a_1\,a_0\,b_2\,b_1 +
a_1\,a_0\,b_1\,b_0 \\
p_2 &=& a_2\,b_0 + a_1\,b_1 + a_0\,b_2\ + a_1\,a_0\,b_1\,b_0\\
p_1 &=& a_1\,b_0 + a_0\,b_1\\
p_0 &=& a_0\,b_0
\end{eqnarray*}

Again let~$I_0$ be the ideal of vanishing polynomials in~$R$.
In this case, the ideal~$I_0$ is generated by the field equations~$x^2-x=0$
for every variable~$x$. Now we compute a \Groebner basis of $I$ in the ring $R/I_0$.
In this ring every ideal is principal (cf. \thmref{unique-generator}) and hence its reduced \Groebner basis will consist of just
one polynomial. Moreover, $I = \skalar{1}$ if and only if its reduced \Groebner basis is~$\{1\}$ and this is equivalent to
the zero set of all polynomials in~$I$ being empty,
and therefore if and only if the property $P$ holds.

\subsubsection{Modelling advantages and disadvantages}

%
Both \modell{}ing approaches presented in
\secref{sec:wordlevel} and \secref{sec:bitlevel} have 	strengths and
weakenesses.
On the one hand, the word-level formulation of verification  problems  as
polynomial systems over~\Ztwon\ leads to fewer variables and equations.
The equations of arithmetic blocks, like multiplier and adder blocks,
are given in a natural and human readable way.
However, not all formul\ae{} on word-level (for example bitwise~\texttt{and},
\texttt{or}, and~\texttt{exclusive-or})
may be coded by polynomial equations.
Therefore, full strength will need bit-level encoding of some variables.
Another drawback are the coefficients from~\Ztwon, which is a ring with
zero-divisors and not a
field. Hence, one cannot rely on valueable properties of
fields, like the algebraic closure.

Since~$\Ztwo$ is a field, these restrictions do not exists for polynomials over~$\Ztwo$,
which can be used for formulation of arbitrary bit-level equations. Moreover,
since the coefficients are restricted to be one or zero, they  need not to be
stored at all.
Hence, a \specialised data structure is possible, which is tailored to suit
this application task.
On the other hand, contrary to the word-level case, bit-level formulations carry many variables and
equations. The number of them may grow exponentially even for some applications
which can be handled easily over~\Ztwon.

As a result from these considerations, research was done for both
approaches. In the following, we present the different strategies and solutions for both,
the word-level and bit-level approach,  in the appropriate algebraic setting.

\section{Standard bases over rings}\label{sec:rgb}
\subsection{Basic definitions}
In this paragraph we outline the general theory of standard bases for ideals or modules
over a polynomial ring $\cR[x_1,\dots,x_n]$ where $\cR$ is any commutative Noetherian
ring with 1. We do not require that the monomial ordering is a well-ordering, that is we
treat the case of standard bases in the localization of $\cR[x_1,\dots,x_n]$ as well (for
a full treatment cf.~\cite{meinphd}).
Gr\"{o}bner bases over $\cR[x_1,\dots,x_n]$ (\ie the case of well-orderings) have been treated
previously (cf.~\cite{adams, kb-appr}) but never for non well-orderings.
Since we are mainly interested in the case $\cR = \Ztwon$ we allow $\cR$ to have
zero-divisors. Moreover, since we are interested in practical application
to real world formal verification problems, we have to develop the theory
for $\cR = \Z_m$ with special care. The ring $\Z_m$ allows special algorithms which dramatically improves the performance
of Gr\"{o}bner bases computations against generic implementations for general rings.

We recall some algebraic basics, including classical notions for
the treatment of polynomial systems, as well as basic definitions and results
from computational algebra. For an exhaustive textbook  about the
subject, when the ground ring $\cR$ is a field, we refer to~\cite{singintro} and the references therein.

Let~$\cR[\monom] = \cR[x_1,\dots,x_n]$ be the polynomial ring over~$\cR$, equipped
with an arbitrary monomial ordering~$<$, \ie global (well-ordering), local or mixed (cf.~\cite{singintro}).
Further $\cR[\vars]_<$ denotes the localization
of $\cR[\vars]$ by the multiplicatively closed set
$$S_< = \{f\in\cR[\vars]\wo\{0\}~|~\LM{f} = 1 \wedge \LC{f}\in \cR^*\},$$
where $\cR^*$ is the group of units of $\cR$ and $\lm$ respectively $\lc$ denote the leading monomial
respectively the leading coefficient w.r.t. $<$, as defined in~\cite{singintro}.
Then $$\wR:=\cR[\vars]_<=\left\{\left.\frac{f}{g}~\right|~f\in\cR[\vars], g\in S_<\right\}.$$
%
%

Also, consider a partition of the ring variables~$\{\vars, \mathbf{y}\}=\{x_1,\dots,x_n,y_1,\dots,y_m\}$.
A monomial ordering over~$\cR[\vars, \mathbf{y}]$ is called an \defemph{elimination
ordering} for $\vars$,
if $x_i>t$ for each~$i$ and for every monomial~$t$ in~$\cR[\mathbf{y}]$.
\pagebreak

\begin{defn}\index{leading!term}\index{leading!coefficient}\index{leading!monomial}\index{leading!ideal}
Let $I\subset\wR=\cR[\vars]_<$ be an ideal and $f$ an element in $\wR$. Choose $u\in S_<$ such that~$\LC{u}=1$
and~$u\cdot f$ is a polynomial $a_0\cdot \monom^{\alpha_0}+\dots+a_n\cdot \monom^{\alpha_n}\in \cR[\vars]$
\label{def:basic_defs}
with $a_0\neq0$ and $x^{\alpha_0} > x^{\alpha_i}$ for all $i\neq 0$ with $a_i\neq 0$ (which is always possible).
Then we define
\begin{align*}
\LT{f} &= a_0\cdot \monom^{\alpha_0}&&\text{leading term of $f$}\\
\LM{f} &= \monom^{\alpha_0}&&\text{leading monomial of $f$}\\
\LC{f} &= a_0&&\text{leading coefficient of $f$}\\
\LE{f} &= \alpha_0&&\text{leading exponent of $f$}\\
\LI{I} &= \skalar{\LT{f}~|~f\in I}_{\cR[\vars]}&&\text{leading ideal of $I$}\\
\LMI{I} &= \skalar{\LM{f}~|~f\in I}_{\cR[\vars]}&&\text{leading monomials ideal of $I$}\\
\V(I) &= \{\vars~|~\forall f\in I: f(\monom) = 0\} &&\text{common zeroes or variety of $I$}\\
\I(V) &= \{f ~|~\forall \vars \in V: f(\monom) = 0\} &&\text{vanishing ideal of $V\subset\cR^n$}\\
\support(f) &= \{\monom^{\alpha_i}~|~a_i\neq 0\}&&\text{support of $f$}\\
\tail(f) &= f - \LT{f} &&\text{tail of $f$}
\end{align*}
If the monomial order~$<$ is global then~$u=1$. If~$<$ is not global the leading coefficients and
the leading terms are well defined, independent of the choice of $u$.
\end{defn}

\begin{defn}\index{grobner basis@\Groebner basis}\index{grobner basis@\Groebner basis!strong}\index{standard basis@standard basis}
\index{standard basis@standard basis!strong}\label{defn_std_basis}
Let $I\subset\wR=\cR[\vars]_<$ be an ideal. A finite set $G\subset R$ is called a
\defemph{standard basis} of $I$ if
\[G\subset I\text{ and }\LI{I}=\LI{G}.\]
That is, $G$ is a standard basis, if the leading terms of $G$
generate the leading ideal of~$I$. $G$ is called a
\defemph{strong standard basis} if, for any $f\in I\wo\{0\}$, there
exists a $g\in G$ satisfying $\LT{g}|\LT{f}$. If $<$ is global we will call standard bases also \defemph{\Groebner bases}.
A finite set $G\subset \wR$ is called standard resp. \Groebner basis, if $G$ is a standard resp. \Groebner basis of $\skalar{G}_\wR$,
the ideal generated by~$G$.
\end{defn}

\begin{rem}
If $\cR$ is a field, than $\LI{I} = \LMI{I}$, but due to non-%
invertible coefficients, in general only
$\LI{I}\subset\LMI{I}$ holds.
\end{rem}

Next, the notion of $t$-representations is introduced, as
formulated in~\citep{Bec93}.
While this notion is mostly equivalent to using syzygies, it helps to understand the
correctness of the algorithms.
\begin{defn}[$t$-representation]
\label{defn:t-repr}
Let $t$ be a monomial and consider elements $$f, g_1, \ldots, g_m, h_1, \ldots , h_m\in \cR[\vars]_<=\wR$$ with
$f = \sum_{i=1}^{m}h_i \cdot g_i$.
Then the sum is called a \defemph{$t$-representation} of $f$ with respect to $g_1, \ldots ,g_m$
if
$$\lm(h_i \cdot g_i) \leq t \text{ for all }i\text{ with } h_i\cdot g_i\neq0\mfstop$$
\end{defn}
\begin{exmp}
Let the monomials of $\cR[x,y]$ be lexicographically  ordered
($x > y$) and
$g_1= x^2, g_2=x^5-y, f=y$.
Then $f=x^3 g_1 - g_2$ is a $x^5y^5$-representation of~$f$.
\end{exmp}

\begin{notation}
Given a representation $p=\sum_{i=1}^{m}h_i\cdot f_i$ with
respect to~$f_1, \ldots f_m$, we may shortly say that $p$ has a
\defemph{nontrivial $t$-representation}, if a $t$-representation of $p$ exists with
$$t<\max\{\lm(h_i\cdot f_i)\vert h_i \cdot f_i\neq 0\}.$$
Note that there exists no $t$-representations with $t<\lm(p)$.
Further, we say that an arbitrary $g$ has a \defemph{standard representation} with respect to~$\{f_i\}$, if
it has a $\LM{g}$-representation.
\end{notation}

\subsection{Normal forms}\label{sec:nf}

{
\renewcommand{\theenumi}{\roman{enumi}}
\renewcommand{\theenumii}{\arabic{enumii}}
\renewcommand{\labelenumii}{(\theenumii)}

\begin{defn}
Let $\calG$ be the set of all finite subsets $G$ of $\wR = \cR[\vars]_<$. A map
\[\nf:\wR\times\calG\to\wR, (f,G)\mapsto\NF{f}{G}\]
\begin{enumerate}
\item is called a \defemph{normal form}\index{normal form} on $\wR$ if, for
all $G\in\calG$,
\begin{enumerate}
\setcounter{enumii}{-1}
\item $\NF{0}{G}=0$,
\end{enumerate}
and, for all $f\in\wR$ and $G\in\calG$,
\begin{enumerate}
\item $\NF{f}{G}\neq 0\Rightarrow\LT{\NF{f}{G}}\not\in\LI{G}$ and
\item $r:=f-\NF{f}{G}$ has a
standard representation with respect to $G$.
\end{enumerate}
\item is called a \defemph{weak
normal form},\index{normal form!weak} if instead of $r$ we just
require that the polynomial~$r'=uf-\NF{f}{G}$ for a unit $u\in\units{\wR}$ has a
standard representation with respect to $G$.

\item is called
\defemph{polynomial weak normal form}\index{normal form!polynomial} if it is a weak normal form and whenever
$f\in\cR[\vars]$ and $G\subset\cR[\vars]$, there exists a unit
$u\in\units{\wR}\cap\cR[\vars]$, such that $uf-\NF{f}{G}$ has a standard
 representation $\sum_{i=1}^n a_ig_i$ w.r.t. $G=\{g_1,\dots,g_n\}$ with $a_i\in\cR[\vars]$.
\end{enumerate}

\end{defn}
}


\begin{rem}
Polynomial weak normal forms exists for arbitrary Noetherian rings and are computable if linear equations
over $\cR$ are solvable (\thmref{thm_NFexist}).
\end{rem}

\begin{defn}
We call a normal form~$\NF{\cdot}{\cdot}$ \defemph{reduced},
if for all~$f\in R$ and~$G\in\calG$ the leading terms of elements from~$G$
do not divide any term of~$\NF{f}{G}$.
Further we call~$G$ a \defemph{reduced \Groebner basis},
if no term from~$\tail(g)$ for any~$g\in G$ is divisible by a leading term of an element of $G$.
\end{defn}

{
\renewcommand{\spoly}[1]{\sPoly\left( #1\right)}         
\renewcommand{\crit}[1]{\Crit\left( #1\right)}           
\renewcommand{\e}{\mathbf{e}}         
Now we introduce an algorithm for computing a polynomial weak normal form for any monomial ordering,
given we are able to solve an arbitrary linear equation in the coefficient ring $\cR$. To ensure correctness
and termination we need to introduce the concept of the $\ecart$ of a polynomial.

\begin{defn}\index{ecart@$\ecart$}
Let $f\in \wR\wo\{0\}$ be a polynomial. Then the \defemph{ecart} is defined by
\[\ecart f = \deg f - \deg \LM{f}.\]
\end{defn}

We introduce a monomial order $<_h$ on $\cR[t, \vars]$ where $t$ is a new variable via
\begin{align*}
t^p\monom^\alpha <_h t^q\monom^\beta :\Longleftrightarrow &\hspace{8.4pt} p + \card{\alpha} < q + \card{\beta} \text{ or }
\\ & \left(p + \card{\alpha} = q + \card{\beta} \text{ and } \monom^\alpha < \monom^\beta\right).
\end{align*}

This is a well-ordering as there are only finitely many monomials with a given total degree. 


\begin{algorithm}
\caption{Calculating a normal form over coefficient rings}
\label{NFBuch}
\begin{algorithmic}
\REQUIRE $f\in \wR$ a polynomial, $G\subset\wR$ finite, $>$ a monomial ordering
\ENSURE A normal form of $f$
\STATE $T := G$ 
\WHILE{$f \neq 0$ and $\LT{f} \in \LI{T}$} 
\STATE   solve $\LT{f} = \sum\limits_{i=1}^s c_i\, \monom^{\alpha_i}\LT{g_i}$
     \begin{tabular}{l}
         with $\monom^{\beta_i}\LM{g_i}=\LM{f}$,
         \\ $g_i\in T$ and $\max\{\ecart g_i\}$ minimal
    \end{tabular}
  \IF{$\max\{\ecart g_i\} > 0$}
    \STATE $T := T \cup \{f\}$
  \ENDIF
  \STATE $f := f - \sum\limits_{i=1}^s c_i\, \monom^{\beta_i}g_i$
\ENDWHILE
\RETURN $f$
\end{algorithmic}
\end{algorithm}

\begin{thm}\label{thm_NFexist}
The \algref{NFBuch} terminates and computes a norm form, if we can solve linear equation in the coefficient ring $\cR$.
\end{thm}

\begin{rem}
In many cases it is not necessary to solve linear equations during the
normal form computation. These include coefficient
fields~(the classical case), weak 1-factorial rings or principal
ideal domains. The latter case was already treated in  \cite{adams}. Further cases can also be computed without
solving linear equations if we require $G$ to be a
strong \Groebner basis.
\end{rem}

\subsubsection{Weak factorial rings}

In rings with zero-divisors we have in general no decomposition into irreducible
elements. For example in $\Z_{12}$ we have $6 = 3\cdot 6 = 3\cdot3\cdot6 =\dots$. Therefore
the concept of factoriality does not make sense. But there exists a notion of weak factorial
rings where every element can be written as $a=n\cdot a_1^{r_1}\cdot\dots\cdot a_k^{r_k}, r_i\geq 0$ ($n$ not
necessarily a unit), such that $a\,|\,b=m\cdot a_1^{s_1}\cdot\dots\cdot a_k^{s_k}$ iff $r_i\leq s_i$. This will be formalized
below.

Let $\cR$ be a commutative Noetherian ring with $1$ and $\units{\cR}$ the group of units. Denote further by
$\zeros{\cR} = \{a\in \cR~|~\exists b\neq 0~:~a\cdot b = 0\}$, the zero-divisors and by
$\nonunits{\cR} = \cR\wo\units{\cR}$ the non-units in $\cR$.

\begin{defn}
An \defemph{element factorization} $(\nu, P)$ or just $\nu$ for a ring $\cR$ consists of a subset
$P\subset\nonunits{\cR}$ and a map $\nu=(\nu_p)_{p\in P}:\cR \to \N^P$, $\nu_p:\cR\to\N$,
such that for all $a\in\cR$ there exists an element $n\in\cR$ with
\[a = n \cdot \prod_{p\in P}p^{\nu_p(a)} =: n\cdot\primes^{\nu(a)}\]
and $\nu_p(a)\neq 0$ only for finitely many $p\in P$.

A ring $\cR$ with an element factorization $\nu$ is called \defemph{$P$-weak factorial} or just \defemph{weak factorial} if, for all $a,b\in\cR$
\[a~|~b \Longleftrightarrow \nu(a) \leq \nu(b).\]
That is, divisibility in $\cR$ is given by the natural order relation of $\N^P$.
If we want to \emphasise the number
of elements in $P$ (elements in $P$ are also called ``primes''), we say weak $\card{P}$-factorial ring where
$\card{P}$ is the cardinality of $P$.
\end{defn}

\begin{example}
\begin{enumerate}
\item If $\cR$ is a factorial domain and $P$ the set of irreducible elements then $\cR$ is
$P$-weak factorial.
\item The ring of integers modulo a power of a prime number $p$ is a weak
1-factorial ring with $P = \{p\}$.
\item The ring $\Z_m$ is weak factorial with $P = \{p\in\mathbb{P}~|~p~|~m\}$, where $\mathbb{P}$ denotes the set of prime numbers.
\item The ring $\Z$ is a weak $\infty$-factorial ring with $P=\mathbb{P}$ and $\nu=\nu^\Z$ the map which
associate to $a\in\Z$ the exponents of the prime decomposition of $a$.
\item The ring $\K[[x]]$, $\K$ a field, is weak factorial with $P=\{x\}$.
\end{enumerate}
\end{example}

\begin{rem}\label{rem:nudef}
For the case of $\Z_m$ with $m=p_1^{e_1}\cdots p_n^{e_n}$, we define $\nu$ as
\[\nu_{p_i}(\underline{a}):=\nu_i(a) = \min\{\nu^\Z_{p_i}(a), e_i\}\]
where $a\in\Z$ represents $\underline{a}\in\Z_m$. E.g. in $\Z_{12}$ we have $12=2^2\cdot 3^1$
and therefore $\nu_3(9)=1$ and $9=3\cdot 3=n\cdot 3^1$. Further in this case~$\nu$ has
the following properties:
\end{rem}

\begin{proposition}
Let $\nu$ be defined for $\Z_m$ as in \remref{rem:nudef}. Then we have
\begin{enumerate}
\item $\nu$ is well-defined, that is $\nu(a) = \nu(a + k\cdot m)$ for all $a,k,m\in\Z$.
\item $\nu$ is saturated multiplicative, that is $\nu_i(a\cdot b)=\min\{\nu_i(a)\cdot\nu_i(b), e_i\}$,
\item $\nu_i(a+b) = 0$ if $\nu_i(a) > 0$ and $\nu_i(b) = 0$,
\item $\nu(a) = 0 \Leftrightarrow \underline{a}\in\units{\Z_m}$ and
\item $\nu$ is \defemph{nice weak factorial}, that is,
$\forall \underline{a}\in\Z_m~ \exists \underline{u}\in\units{\Z_m}\,:\,\underline{a}=\underline{u}\cdot\primes^{\nu(\underline{a})}$.
\end{enumerate}
\end{proposition}
\begin{proof}
The first four properties follow easily from the valuation properties of $\Z$ with $\nu^\Z$. For the
last one let $\underline{a} = \underline{n} \cdot \primes^{\nu(\underline{a})}$.
At first notice, that $\nu_{p_i}(\underline{n}) > 0$ is only possible, if $\nu_{p_i}(\underline{a}) = e_i$. Hence consider
\[u = n + \frac{m}{\primes^{\nu(a)}}\cdot\mathop{\prod_{e_i > 0}}_{p_i \nmid\, n}p_i.\]
Now $\nu(\underline{u}) = 0$ and therefore $\underline{u}\in\units{\Z_m}$. Further
$\underline{u}\cdot\primes^{\nu(\underline{a})}=\underline{a}$.
\end{proof}

\begin{rem}
One can show that in our definition the elements of $P$ are irreducible and that $\cR$ is a weak
unique factorization ring (UFR) in the sense of Agarg\"{u}n \cite{agar1999} and therefore
a \generalisation of the notions from Bouvier-Galovich \cite{galovich1978, bouvier} and
Fletcher \cite{fletcher1969} (cf. \cite{agar1999}). Nevertheless we prefer our definition, as it emphasis the divisibility relation.
\end{rem}

\begin{rem}
If $\cR$ is a principal ideal ring, then it is isomorphic to a finite product~\cite{samzar}
of principal ideal domains,
hence factorial domains, and finite-chain rings (cf. \cite{samzar}), which are weak 1-factorial. Therefore we can
compute \Groebner basis in polynomials rings over the factors
and lift them to $\cR[\vars]$. This is described in the work of G.\ Norton and A.\ Salagean \cite{norton}. Below
we show that computation in the
ring itself is feasible.
\end{rem}

\begin{defn}\index{greatest common devisor}\index{least common multiple}\index{gcd}\index{lcm}
Let $\cR$ be a weak factorial ring and $a_1, \dots, a_n\in \cR$. Then we define~(with $\max, \min$ component-wise)
\begin{align*}
\ggt{a_1,\dots,a_n} &= \primes^{\min\{\nu(a_1),\dots,\nu(a_n)\}}\quad\text{and}\\
\kgv{a_1,\dots,a_n} &= \primes^{\max\{\nu(a_1),\dots,\nu(a_n)\}}.
\end{align*}
\end{defn}

\begin{rem}
This definition of $\gcd$ and $\lcm$ fulfills the universal properties of the greatest common divisor
and the least common multiple. But notice that, for arbitrary rings,
the $\gcd$ and $\lcm$ are not unique up to units. However, in the case of $\Z_m$ this holds:
\[a|b\wedge b|a \Rightarrow \exists u\in\units{\Z_m}~:~a=u\cdot b.\]
\end{rem}

\begin{lem}
Let $\cR$ be a weak factorial principal ring. Then
\begin{align*}
\skalar{a_1, \dots, a_n} &= \skalar{\ggt{a_1,\dots,a_n}},\\
\skalar{a_1} \cap \dots \cap \skalar{a_n} &= \skalar{\kgv{a_1,\dots,a_n}}.
\end{align*}
\end{lem}

\begin{proof}
Follows directly from the definition of weak factorial and $\gcd$, respectively $\lcm$, and their universal properties.
\end{proof}

\begin{lem}\label{lem_CND_eq}
Let $\cR$ be a weak 1-factorial principal ring with prime~$\eta$ and let~$c, a_1, \dots,
a_s\in\cR\wo\{0\}$. Then the following are equivalent.
\begin{itemize}
\item The equation $c = a_1x_1+\dots+a_sx_s$ is solvable.
\item There exists an $j\in\{1,\dots, s\}$ and  $x\in\cR$, such
that $c=a_jx$, \ie $a_j|c$.
\end{itemize}
\end{lem}

\begin{proof}
The first statement is equivalent to
\begin{align*}
                & c\in\skalar{a_1,\dots, a_n} \\
\Leftrightarrow & \ggt{a_1,\dots,a_n}~|~c \\
\Leftrightarrow & \min\left\{\nu(a_1),\dots,\nu(a_n)\right\} \leq \nu(c) \\
\Leftrightarrow & \exists a_i~:~\nu(a_i) \leq \nu(c), \text{ as } \Image{\nu}\subset\N \\
\Leftrightarrow & c\in\skalar{a_i}
\end{align*}
which is equivalent to the second statement.
\end{proof}

\begin{cor}
Let $\cR$ be a weak 1-factorial principal ring.
Then, solving linear equations over $\cR$ can be reduced to tests for
divisibility. Moreover, every standard
basis over $\cR[\vars]_<$ is a strong standard basis.
\end{cor}



\subsection{Computing standard bases}\label{sec:comp_std_bases}

Let $\cR$ be a commutative Noetherian ring with $1$.

\begin{defn}\index{syzygy module}
Let $R$ be a ring and $A\in \wR^{s\times t}$ a matrix considered as a linear map $\wR^s\to\wR^t$. The kernel of
$A$ is a submodule of $\wR^s$. It is called the \defemph{syzygy module} of
$A$. If $A=(f_1, f_2, \dots, f_s)\in\wR^{s\times 1}$, then
\[\Syz{f_1, \dots, f_s}=\ker(A)=\{(h_1, \dots, h_s)\in\wR^s ~|~\sum h_if_i=0\}.\]
\end{defn}

\begin{thm}(Buchberger's criterion)\label{thm_buchcrit}\index{Buchberger's criterion}
Let $I\subset \wR=\cR[\vars]_<$ be an ideal and $G=\{g_1,\dots,g_s\}\subset I$.
Further let $\NF{-}{G}$ be a weak normal form on $R$ with
respect to $G$. Then the following statements are equivalent:
\begin{enumerate}
\item $G$ is a standard basis of $I$.
\item $\NF{f}{G}=0$ for all $f\in I$.
\item Each $f\in I$ has a standard representation with respect to
$G$.
\item $G$ generates $I$ and for every element $h$ with
\[h\in{\Syz{\LT{g_i}|i=1,\dots,s}},\]
$\NF{h_1g_1+\dots+h_sg_s}{G}=0$.
\end{enumerate}
\end{thm}

\begin{proof}
The implications $1\Rightarrow2\Rightarrow3\Rightarrow4\Rightarrow1$
can be shown as in the classic case. The classical proof
can be found either in \cite{singintro} (general orderings)
or in \cite{adams} (global orderings).
\end{proof}

%


To specialize further for the case of weak factorial principal rings we modify the
classical notion of an $s$-polynomial.

\begin{defn}\index{s-polynomial@$s$-polynomial}
Let $f,g\in \wR\wo\{0\}$. We define the \defemph{$s$-polynomial} of $f$
and $g$ to be
\[\spoly{f,g}:=\frac{\kgv{\LT{f},\LT{g}}}{\LT{f}}f-\frac{\kgv{\LT{f},\LT{g}}}{\LT{g}}g.\]
\end{defn}

\begin{rem}
This definition is not equivalent to
\begin{multline*}
\spolyr{f,g}=\\\LC{g}\frac{\kgv{\LM{f},\LM{g}}}{\LM{f}}f-\LC{f}\frac{\kgv{\LM{f},\LM{g}}}{\LM{g}}g.
\end{multline*}
For example let $f=2\,x-2\,y, g=2\,y-z$ in $\Z_4[x,y,z]$. Then we
get \linebreak$\spolyr{f,g}=x\,z\neq-2\,y+z\,x=\spoly{f,g}$. That is, we can
loose terms just by multiplying with a constant, \eg if $2\,x+y\in I$ for some ideal $I$,
then~$2\,y\in\LI{I}$. Therefore we have to look for further generators
of the syzygies, the classical $s$-polynomials are not sufficient.
\end{rem}

\begin{defn}
Let $\cR$ be a principal ring and $a\in\cR$. The annihilator of $a$, $\text{Ann}(a)=\{n\in\cR~|~a\cdot n = 0\}$
is an ideal in $\cR$ and is hence generated by one element, which we denote by $\zeroDiv{a}$.
\end{defn}

Due to zero-divisors we define the $s$-polynomial also for pairs $(f,g)$ with one component being $0$.

\begin{defn}\index{s-polynomial!extended}
Let $f\in\wR\wo\{0\}$. We define the \defemph{extended s-polynomial}
of $f$ to be
\[\spoly{0,f}=\spoly{f,0}:= \zeroDiv{\LC{f}}\cdot f.\]
\end{defn}

\begin{algorithm}
\caption{Computes a standard basis of $I$} \label{Std_CNU}
\begin{algorithmic}
\REQUIRE\
\\$I$ a finite set of polynomials,
\\$>$ a monomial ordering, $\nf$ a weak normal form
\ENSURE $G$ is a standard basis of $I$
\STATE $G := I$
\STATE $P := \{(f,g)~|~f,g\in S, f\neq g\} \cup \{(0,f)~|~f \in G\}$, the pair set
\WHILE{$P \neq \emptyset$}
  \STATE choose $(f,g)\in P$
  \STATE $P := P\wo\{(f,g)\}$
  \STATE $h := \NF{\spoly{f,g}}{G}$
  \IF{$h\neq 0$}
    \STATE $P := P\cup\{(h,f)\,|\,f\in G\} \cup \{(0,h)\}$
    \STATE $G := G\cup\{h\}$
  \ENDIF
\ENDWHILE
\RETURN $G$
\end{algorithmic}
\end{algorithm}

\subsubsection{Buchberger's criterion and the syzygy theorem}

In the following we assume $\cR$ to be a weak factorial principal ring.
Termination of \algref{Std_CNU} is an easy consequence of the Noetherian property of the ring $\wR$.
To present the theorem, which implies the correctness of \algref{Std_CNU} we need to
introduce some terminology. We fix a set of generators $G=\{f_0,f_1,\dots,f_k\}$ of an
ideal $I$ with $f_0 = 0$.

First assume that a set
$J\subset\{(i,j)~|~0\leq j < i\leq k\}$ is given with
$$\NF{\spoly{f_i,f_j}}{G} = 0 \text{ for }(i,j)\in J.$$
For $0\leq i<j\leq k$ let $\LT{f_i} = c_i \monom^\alpha_i$ and define:
\begin{align*}
m_{ji} &= \frac{\kgv{c_i,c_j}}{c_i}\cdot \frac{\kgv{\monom^{\alpha_i},\monom^{\alpha_j}}}{\monom^{\alpha_i}} 
        = \frac{\kgv{\LT{f_i},\LT{f_j}}}{\LT{f_i}}\\
m_{0i} &= \zeroDiv{c_i}\\
\spoly{f_i,f_j} &= m_{ji} f_i - m_{ij} f_j \\
\spoly{f_i,f_0} &= m_{0i} f_i \text{ as } f_0 = 0 \text{ (set also } m_{i0} = 0\text{)}\\
\spoly{f_i,f_j} &= \sum_{\nu=1}^k a_\nu^{(ij)} f_\nu\text{ the standard representation for }(i,j)\in J\\
s_{ij} &= m_{ji} \e_i - m_{ij} \e_j - \sum_{\nu=1}^k a_\nu^{(ij)} \e_\nu\in \Syz{I}\text{ for }(i,j)\in J
\end{align*}

The elements $m_{0i}$ and $s_{i0}$ correspond to the new $s$-polynomials, which occur due to zero divisors.

\begin{thm}[Buchberger's criterion]\label{thm_buchsyz}\index{Buchberger's criterion}\index{Syzygy theorem}
Let $G=\{f_0, f_1,\dots,f_k\}$ be a set of generators of $I\subset \wR$ with $f_0 = 0$. Further let $J\subset\{(i,j)~|~0\leq i < j \leq k\}$ be such that $\skalar{{m_{ij}\e_j~|~(i,j)\in J}} = \skalar{{m_{ij}\e_j~|~0\leq j < i \leq k}}$.
If
$$\NF{\spoly{f_i, f_j}}{G_{ij}} = 0\text{ for }(i,j)\in J$$ and some $G_{ij}\subset G$ then
\begin{enumerate}
\item[(a)] $G$ is a standard basis of $I$ (Buchberger's criterion) and
\item[(b)] $S := \{s_{ij}~|~(i,j)\in J\}$ generates $\Syz{I}$.
\end{enumerate}
\end{thm}

For a proof we refer to \cite{meinphd}.

\begin{remark}
The set $S$ is a standard basis of $\Syz{I}$ with respect to the
Schreyer ordering (definition of the Schreyer ordering cf. \cite{singintro}).
\end{remark}

\begin{cor}
\algref{Std_CNU} terminates and is correct.
\end{cor}

\begin{remark}
If $f$ and $I$ are polynomial and if $\nf$ is a polynomial weak normal form in \algref{Std_CNU}
than $G$ is a standard basis of $\skalar{I}_R$ consisting of polynomials.
\end{remark}
%
%
Also, the $t$-representations of \defref{defn:t-repr} can be \utilised for a standard basis test as given below.
%
\begin{theorem}
Let $F=(0, f_1, \ldots, f_k)$, $f_i \in \cR[\vars]$, be a polynomial system.
If~$\spoly{f,g}$ has a nontrivial $t$-representation \wrt~$F$
for each~$f, g\in F$, then~$F$ is a \Groebner basis.
\end{theorem}

\begin{proof}
The theorem can be proved similar as in~\citep{Bec93}. A more sophisticated
version of this theorem can be formulated and proven likewise
to~\citep[p. 142]{singintro}.
\end{proof}



\subsubsection{Criteria for $s$-polynomials}

In order to compute non-trivial standard bases in practise, we like to have criteria to omit
unnecessary critical pairs. This improves the time and space requirement
of the Buchberger algorithm as in the classical case.

\begin{lem}[Product criterion]
Let $f,g\in\wR=\cR[\vars]_<$ with $\LM{f}$ and $\LM{g}$ relatively prime.
Further let $\LC{f}$ and $\LC{g}$ be a unit, then
\[\NF{\spoly{f,g}}{\{f, g\}}=0.\]
\end{lem}
\begin{proof}
No change of the classical proof is needed. However, the strong product criterion,
which gives an if and only if statement,
is not extendable to the general case.
\end{proof}


\begin{example}
The polynomials $4x+y$ and $y^2+2\,z\in\Z_8[x,y,t]$ will reduce to zero by a
sharper product criterion (not given here). In
contrast $4\,y+x^3+1$ and~$x^5+2\,x^2$ will reduce to $2\,x^2$, which is
not reducible by either of the polynomials nor their extended
$s$-polynomials.
\end{example}

\begin{lem}[Chain criterion]\label{lem_chain_theo}
With the notations of \thmref{thm_buchsyz} let
$\LT{f_i} = c_i\,\monom^{\alpha_i}$,
$\LT{f_j} = c_j\,\monom^{\alpha_j}$, and $\LT{f_l} = c_l\,\monom^{\alpha_l}$ with $i>j>l$. If $c_j\,\monom^{\alpha_j}$ divides
$\kgv{c_i\,\monom^{\alpha_i},c_l\,\monom^{\alpha_l}}$ then $m_{li}\,\e_i\in\skalar{m_{ji}\,\e_i}$. In particular, if $s_{ij}, s_{jl}\in S$
then $S\wo\{s_{il}\}$ is already a standard basis of $\Syz{I}$ and $S\wo\{s_{il}\}$ generates $\Syz{I}$.
\end{lem}

\begin{proof}
The divisibility of $\kgv{c_i\,\monom^{\alpha_i},c_l\,\monom^{\alpha_l}}$ by $c_j\,\monom^{\alpha_j}$ implies
\[\kgv{c_i\,\monom^{\alpha_i},c_j\,\monom^{\alpha_j}}~|~\kgv{c_i\,\monom^{\alpha_i},c_l\,\monom^{\alpha_l}}.\] Dividing
both sides by $c_i\,\monom^{\alpha_i}$ yields $m_{ji}~|~m_{li}$.
\end{proof}

The following criterion is new and quite useful in practise.

\begin{lem}\label{lem_zero_crit}
With the notations of \thmref{thm_buchsyz} let $\LT{f_i} = c_i\,\monom^{\alpha_i}$ and $\LT{f_l} = c_l\,\monom^{\alpha_l}$ with $i>l$.
If $\zeroDiv{c_i}$ divides
$\kgv{c_i, c_l}$ then $m_{li}\,\e_i\in\skalar{m_{0i}\,\e_i}$. In particular, if the special $s_{i0}\in S$ (corresponding to an $s$-polynomial
with one zero entry)
then $S\wo\{s_{il}\}$ is already a standard basis of $\Syz{I}$.
\end{lem}

\begin{proof}
Follows from $m_{0i} = \zeroDiv{c_i}$.
\end{proof} 

}

\section{Boolean \Groebner Basis}\label{sec:bgb}
In the following, we present methods for treating the bit-level formulation of digital
systems as
introduced in \secref{sec:bitlevel}.
First, the notion of Boolean polynomials is given, and a suitable data
structure is motivated.
The next part is addressed to effective algorithms for operations on these
polynomials. Then recent results in the  theory of Boolean \Groebner bases are
presented, including new criteria, which  \minimise the number of critical
pairs. Finally, we sketch a new approach, which improves the algorithms by
exploiting symmetries in the polynomial system.

\subsection{Boolean Polynomials}
\label{sec:poly_as_set}
%
%
\label{sec:def_boole_poly}

In this section we model expressions from propositional logic as
polynomial equations over the finite field with two elements. In this
algebraic language the problem of satisfiability can be approached by
a tailored \Groebner basis computation. We start with the polynomial
ring~$\Ztwoxoneton = \explZtwoxoneton$.

%
Since the considered polynomial functions take only values
from~$\Ztwo$, the condition~$x = x^2$ holds for
all~$x\in\Ztwo$.
%
Hence, it is reasonable to simplify a polynomial in~\Ztwoxoneton \wrt the
\defemph{field equations}
\begin{equation}
\label{eqn:fieldeqns}
x_1^2 = x_1^{\phantom{2}},\, x_2^2 = x_2^{\phantom{2}},\quad\dots\quad,\,
        x_n^2 = x_n^{\phantom{2}}\mfstop
\end{equation}
Let~$\FP=\set{\explfieldequations}$ denote the corresponding set of \defemph{field polynomials}.
The field equations yield a degree bound of one on all variables
occurring in a polynomial in~$\Ztwoxoneton$ modulo~\FP.

\begin{definition}[Boolean Polynomials]
\label{def:boolepoly}
Let~$p\in \Ztwoxoneton$ be a polynomial, \sth
\begin{equation}
\label{eqn:expandedpoly}
p =  a_1 \cdot x_1^{\nu_{11}} \cdot \ldots \cdot x_n^{\nu_{1n}}+ \ldots +
     a_m \cdot x_1^{\nu_{m1}} \cdot \ldots \cdot x_n^{\nu_{mn}}
\end{equation}
with coefficients~$a_i \in \{0, 1\}$. If~$\nu_{ij}\leq 1$ for all~$i,j$,
then~$p$ is called a \defemph{Boolean polynomial}.

The set of all Boolean polynomials in~$\Ztwoxoneton$ is denoted by~\Bool.
\end{definition}


Note that Boolean polynomials can be uniquely identified with a subset of the
power set of~$\{ x_1, \dots, x_n\}$:

\begin{lemma}
\label{lem:poly_as_set}
Let~$R=\Ztwoxoneton$, and~$P=\powerset{x_1, \dots, x_n}$ be the power set of
the set of variables  of~$R$. Then the power set~$\powerset{P}$ of~$P$
is in one-to-one correspondence with
the set of
Boolean polynomials in~$R$ via the mapping $f: \powerset{P} \rightarrow R$
defined
by~$S \mapsto \sum_{s\in S} \left( \prod_{x_\nu \in s}x_\nu \right)$.

\end{lemma}

\begin{pf}
It is obvious, that~$\sum_{s\in S} \left( \prod_{x_\nu \in s}x_\nu \right)\in
\Bool$
for each subset~$S$ of~$P$.
On the other hand,  with the notation of \eqnref{eqn:expandedpoly},
a Boolean polynomial~$p$ is uniquely determined by the fact,
whether a term~$x_1^{\nu_{i1}} \cdot \ldots \cdot x_n^{\nu_{in}}$ occurs in
it, because its coefficents lie in~$\{0, 1\}$.
Moreover, each term is determined by the occurrences of its variables.
Hence, one can assign the set~$S_p= \{s_1, \,\cdots, s_m\}$ to~$p\in\Bool$,
where~$s_k\subseteq\{x_1, \dots, x_n\}$ is the set of variables occurring in the
k-th term of~$p$.
\end{pf}

For practical applications it is reasonable to assume \emph{sparsity}, \ie
the set~$S$ is only a small subset of the power set over the variables.
Even the elements of~$S$ can be considered to be sparse, as usually
only few variables occur in each term.
Consequently, the strategies of the proposed algorithms try to preserve
this kind of sparseness.


The following statements are not difficult to prove, but essential for the whole
theory.

\begin{thm}
    \label{canonical-representatives}
The composition~$\Bool\hookrightarrow  \Ztwoxoneton \twoheadrightarrow \Ztwoxoneton/\fieldideal$ is a bijection.
That is,
the
Boolean polynomials are a canonical system of representatives  of the residue classes
in the quotient ring of~$\Ztwoxoneton$ modulo the ideal of the field
polynomials~$\fieldideal$. Moreover, this bijection provides~\Bool with the
structure of a~$\Ztwo$-algebra.
\end{thm}

\begin{proof}
The map is certainly injective. Since any polynomial can be reduced to a Boolean
polynomial using~\FP, the map is also surjective.
\end{proof}

\begin{defn}
    A function $f:\Ztwo^n\rightarrow\Ztwo$ is called a \defemph{Boolean function}.
\end{defn}

\begin{proposition}
\label{prop:residue_class}
    Polynomials in the same residue class modulo~\fieldideal generate the same
function.
\end{proposition}
\begin{pf}
    Let $p$, $q$ be polynomials with $p-q\in\fieldideal$.
   By \thmref{canonical-representatives} we have
   \[p=b+f_p,q=b+f_q\mcomma\]
  where the first summand $b$ is a common Boolean polynomial and the second
summand lies in $\fieldideal$. The latter evaluates to zero at each point in $\Ztwo^n$.
\end{pf}

\begin{thm}
The map from~\Bool
to the set of Boolean functions~$\{f:\Ztwo^n \rightarrow \Ztwo\}$ by mapping a
polynomial to its polynomial function
is an isomorphism
of $\Ztwo$-vector-spaces. Even more, it is an isomorphism of $\Ztwo$-algebras.
\label{polys-function}
\end{thm}

\begin{pf}
The map is clearly a $\Ztwo$-algebra homomorphism. Injectivity follows from
\thmref{canonical-representatives} together with \propref{prop:residue_class}.
For surjectivity it suffices to see, that both sides have dimension $2^n$.
\end{pf}

\begin{corollary}
    \label{existszeroone}
Every Boolean polynomial $p\neq 1$ has a zero over $\Ztwo$.
Every Boolean polynomial $p\neq 0$ has a one over $\Ztwo$, that is~$p+1$ has a
zero.
\end{corollary}

Recalling \defref{def:basic_defs}, for~$I\subseteq\Ztwoxoneton$ the algebraic
set in~$\Ztwo^n$ defined by~$I$
is denoted by~$\V(I) = \{\vars\in \Ztwo^n~|~\forall f\in I: f(\vars) = 0\}$.

\begin{corollary}
\label{cor:boole1-1}
There is a natural one-to-one correspondence between Boolean polynomials and
algebraic subsets of~$\Ztwo^n$, given by $p\mapsto \V(\pwithfe{p})$.
Moreover, every subset of~$\Ztwo^n$ is algebraic.
\end{corollary}

\begin{proof}
Since~$\Ztwo^n$ is finite, every subset is algebraic.
Let~$\chi_S$ be the characteristic function of a subset~$S \subseteq\Ztwo^n$,
that is~$\chi_S(\vars)=1$ if and only if~$\vars\in S$.
 By  \thmref{polys-function} there is a~$p \in\Bool$ defining~$1+\chi_S$.
Hence, the map is surjective. Moreover, since both sets have the same
cardinality, the results follows.
\end{proof}

After showing the correspondence between Boolean functions and Boolean polynomials we have a look at Boolean formulas, the kind of formulas defining Boolean functions.

\begin{defn}
\label{def:maplogic}
    We define a map $\phi$ from formulas in propositional logic to Boolean
functions,
by providing a translation from the basis system~\texttt{not}~($\neg$), \texttt{or}~($\lor$), \texttt{true}~($\True$).
    For any formulas~$p,q$ we define the following rules
\begin{equation}
\label{eqn:logic2func}
\begin{array}{rcl}
    \phi(p \lor q)&:=&\phi(p)\cdot \phi(q)\\
    \phi(\neg p)&:=&1-\phi(p)\\
    \phi(\True)&:=&0
\end{array}
\end{equation}
    Recursively every formula in propositional logic can be translated into
Boolean functions, as $\{\lor, \neg, \True\}$ forms a basis system in propositional logic.
\end{defn}

\begin{rem}
\begin{enumerate}
\item
It is quite natural to identify $0$ and $\True$ in computer algebra, as we
usually associate to a polynomial $f$ the equation $f=0$, and $f$ being zero is
equivalent to the equation being fulfilled.
\item
For every Boolean function~$f$ there exists a formula~$p$ in propositional
logic, \sth~$\phi(p)=f$. Together with \thmref{polys-function} we obtain that
every formula give rise to a Boolean polynomial, generated by  rules
corresponding to those of \eqnref{eqn:logic2func}.
\end{enumerate}
\end{rem}


We are interested in a representation of Boolean polynomials, whose storage
space scales well with the number of terms and still allows to carry out
vital computations for \Groebner basis computation in reasonable time.
In the next section, a data structure with the desired
properties is presented. Therefore, it can be used to store and
handle the construction of Boolean polynomials proposed in
\lemref{lem:poly_as_set}.

\subsection{Zero-suppressed Binary Decision Diagrams}
\label{sec:bdd}

Binary decision diagrams (BDDs) are widely used in formal verification and model
checking for representing large sets. For instance, they arise
from configurations of Boolean functions and states of
automata which cannot be
constructed efficiently by an enumerative approach.
One of the advantages of BDDs is the performance of basic operations like
intersection and complement. Another major benefit are equality tests, which
can be carried out immediately, as BDDs allow a canonical form.
For a more detailed treatment of the subject
see~\cite{Ghasemzadeh05} and~\cite{berardetal01}.

\begin{definition}[Binary Decision Diagram]

A \defemph{binary decision diagram}~(BDD) is a rooted, directed, and acyclic
graph with two terminal nodes~$\{0, 1\}$ and decision nodes.
The latter have two ascending edges~(high/low or then/else),
each of which corresponding to the assignment of true or false, respectively, to
a given Boolean variable.
In case that the variable order is constant over all
paths, we speak of an \defemph{ordered} BDD.

\end{definition}

This data structure is compact, but easy to describe and implement.
Also, the subset of the power set represented by a BDD
can be recovered easily, by following then- and else-edges.

\begin{definition}
Let~$b$ be a binary decision diagram.
\begin{itemize}
\item
The decision variable associated to the root node of~$b$ is denoted
by~$\mathrm{top}(b)$. Furthermore,
$\thenBranch(b)$ and~$\elseBranch(b)$ indicate
the (sub-)diagrams, linked to then- and else-edge, respectively, of the root
node of~$b$.
\item For two BDDs~$b_1, b_0$, which do not depend on the
decision variable~$x$, the \defemph{if-then-else operator}
$ite(x,b_1,b_0)$ denotes the BDD~$c$, which is obtained by introducing
a new node associated to the variable~$x$, \sth~$\thenBranch(c) = b_1$,
and~$\elseBranch(c) = b_0$.
\end{itemize}
\end{definition}

\begin{figure}[htb!]
\begin{center}
\subfigure[initial diagram]{%
\includegraphics[scale=0.225]{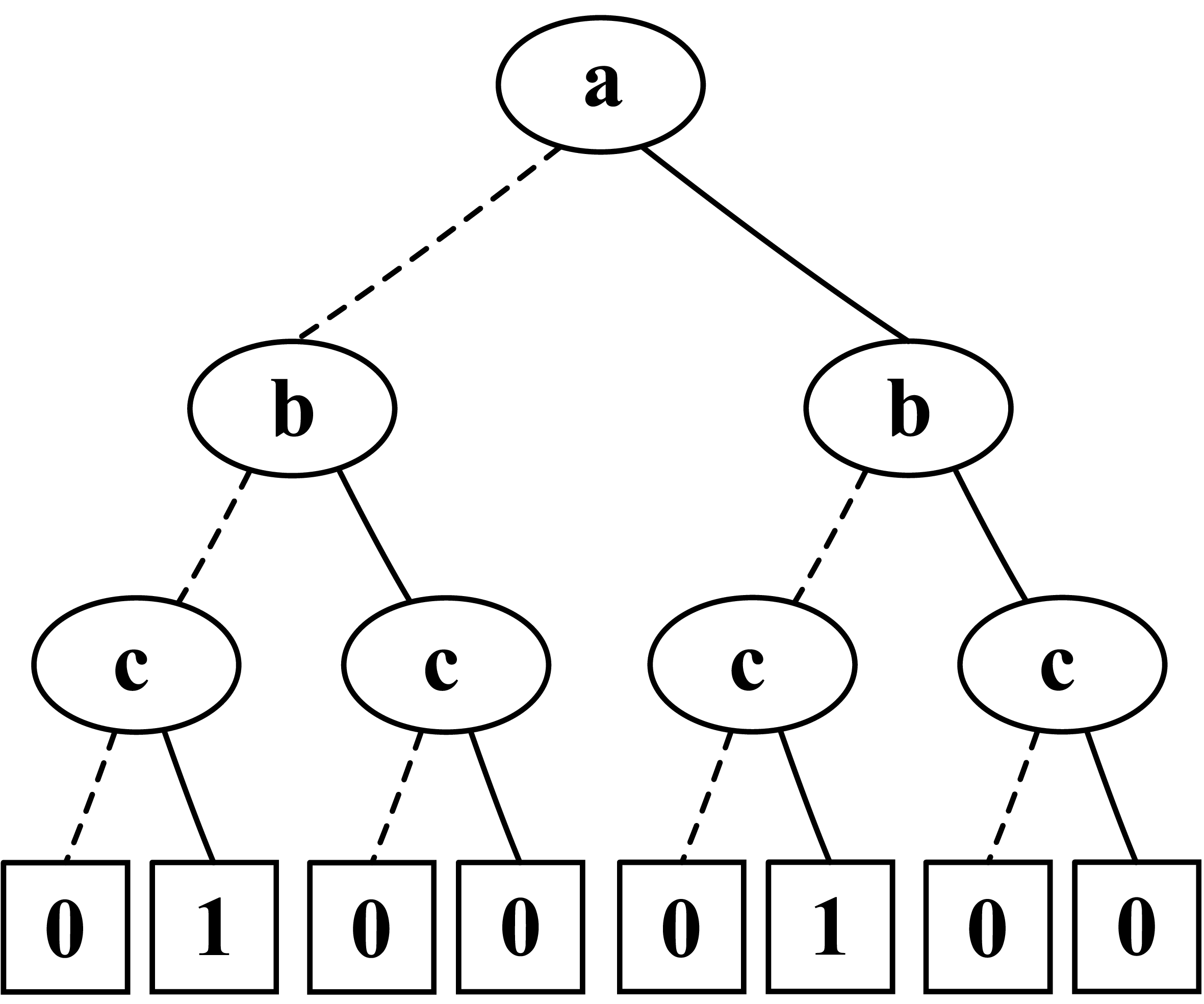}\label{fig:zdd_motiv_1}
}%
\subfigure[subdiagrams merged]{%
\phantom{(b)~}\hspace{1.5em}%
\includegraphics[scale=0.225]{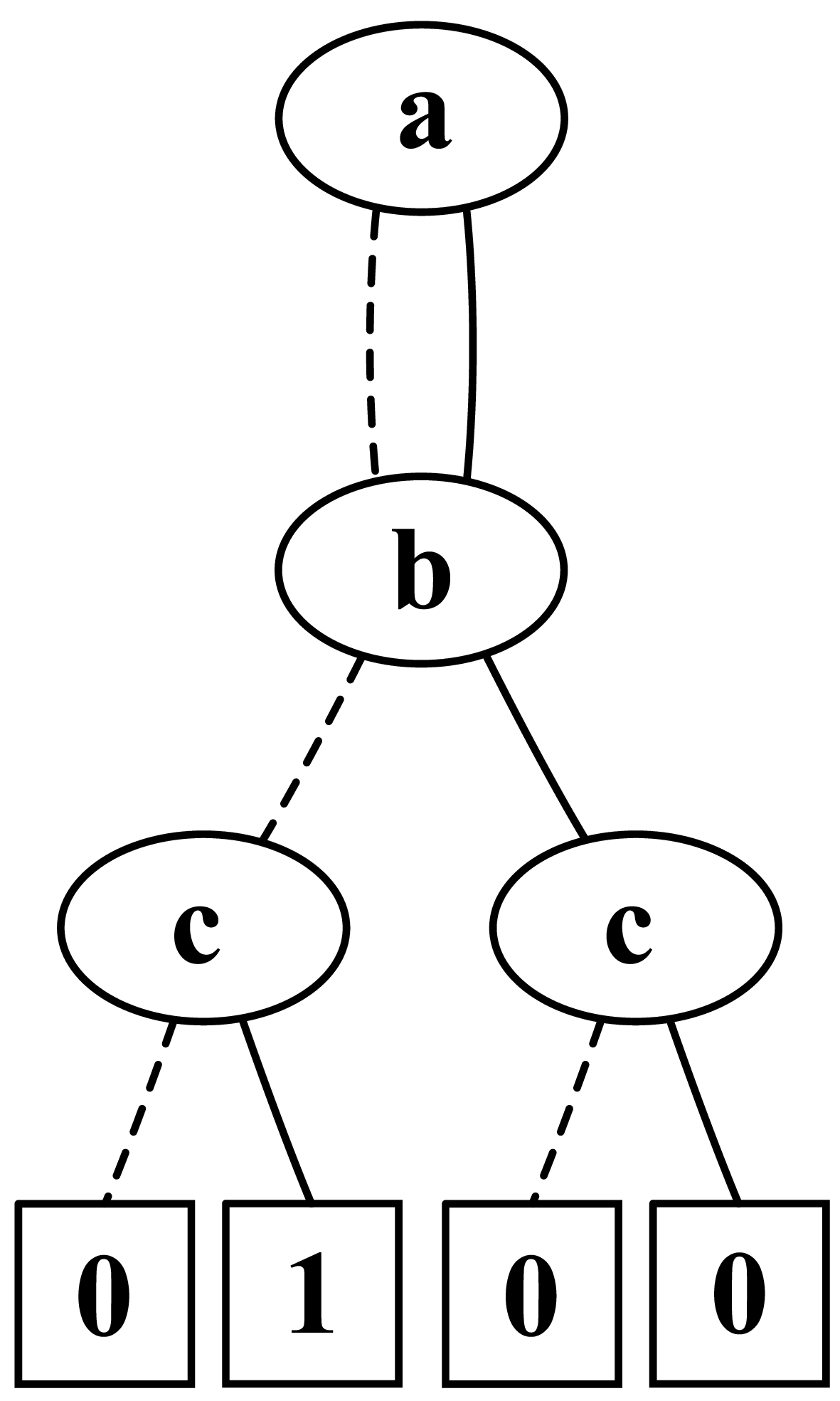}\label{fig:zdd_motiv_2}
\hspace{1.5em}
}%
\subfigure[zero-supressed]{%
\phantom{(c)~}\hspace{1.5em}%
\includegraphics[scale=0.225]{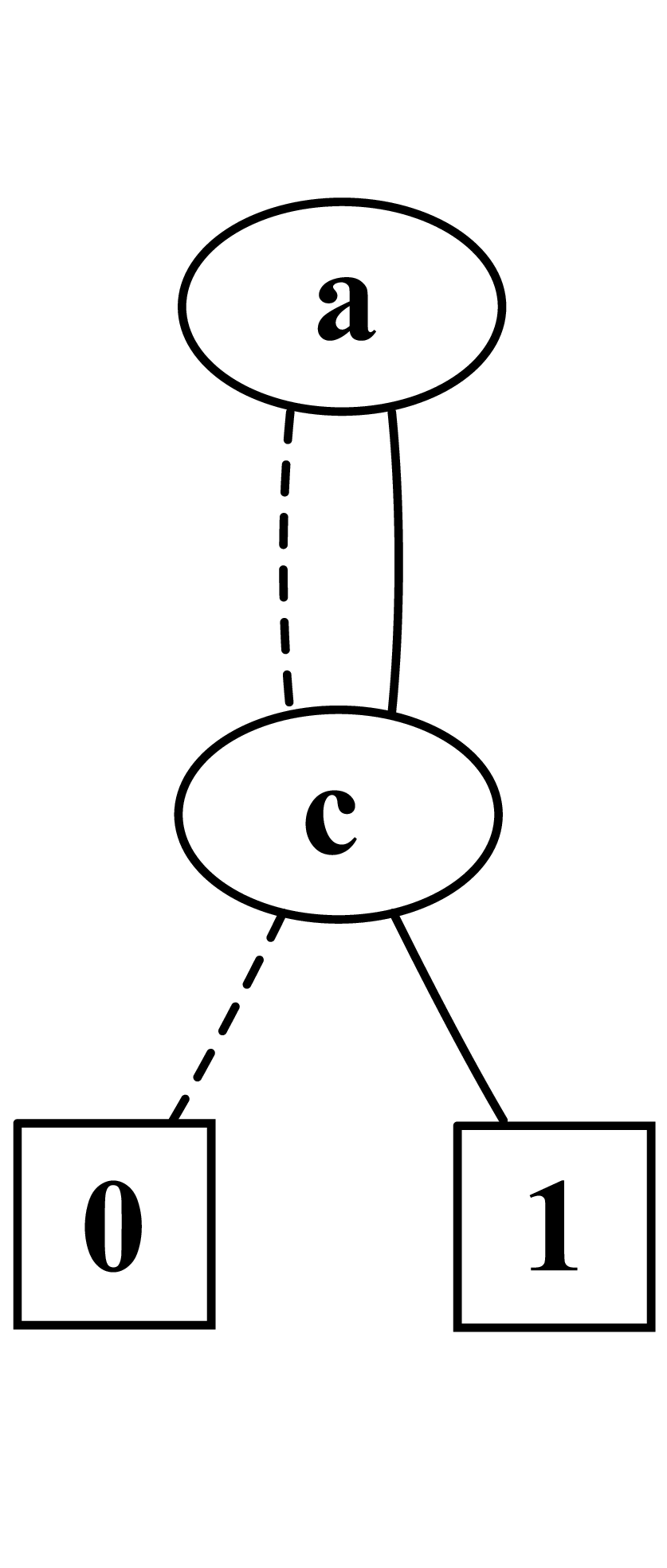}\label{fig:zdd_motiv_3}
\hspace{1.5em}%
}%
\end{center}
\caption{Different kinds of binary decision diagrams representing the polynomial~$a\, c+c$.
Solid/dashed connections marking then/else-edges, respectively.
\label{fig:zdd_motiv}}
\end{figure}

A Boolean polynomial~$p$ can be converted to an ordered BDD using the
following approach. Having variables~$x_1,\ldots,x_n$ the polynomial~$p$  can be
written as~$p= x_1\cdot p_1 + p_0$,
where~$p_1$ and~$p_0$ are Boolean polynomials depending on~$x_2,\dots, x_n$ only.
Therefore,
if we have diagrams~$b_1, b_0$ representing~$p_1$ and~$p_0$, respectively,
the whole diagram is generated by~$ite(x_1, b_1,b_0)$. But~$b_1,
b_0$ can be
obtained by recursive application of the procedure with respect
to~$x_2, \ldots, x_n$. The recursion ends up by a
constant polynomial, which is to be connected to the corresponding terminal
node.
\Figref{fig:zdd_motiv_1} illustrates such a decision diagram  for the
polynomial~$a\,c +c = a\cdot(b\cdot (c\cdot 0 +0)  + (c\cdot 1 +0)  )  +b\cdot(c\cdot 0 +0)  + c\cdot 1 + 0$.
From this example, one can already see, that it is useful to
identify equivalent subdiagrams in such a way that those edges which point to
equal subgraphs are actually linked to the same subdiagram instances.
The merging procedure is sketched in \figref{fig:zdd_motiv_2}.

For efficiency reasons, one may omit variables, which are not necessary
to reconstruct the whole set.
This leads to even more compact representations, which are faster to handle.
A classic variant for this purpose is the
\defemph{reduced-ordered BDD}~(ROBDD, sometimes referred to as ``\emph{the}
BDD''). These are ordered BDDs with equal subdiagrams merged.
Furthermore, a node elimination is applied, if
both descending edges point to the same node.
While the last reduction rule is useful for describing numerous Boolean-valued
vectors, it is gainless for treating sparse sets. For this case,
another variant, namely the ZDD (sometimes also called~ZBDD or~ZOBDD),
has been introduced.

\begin{definition}[ZDD]
Let~$z$ be an ordered binary decision diagram with equal subdiagrams merged.
If those nodes are eliminated whose then-edges point to
the~$0$-terminal, then~$z$ is called a
\defemph{zero-suppressed binary decision diagram}~(ZDD).
\end{definition}

Note, in this case elimination means that a node~$n$ is removed from the
diagram and all edges pointing to it are linked to~$\elseBranch(n)$.
In \figref{fig:zdd_motiv_2} the then-edge of the right node with decision
variable~$c$  is pointing to the 0-terminal. Hence, it can be safely removed,
without losing information.
As a consequence, the then-edge of the $b$-node is now connected to zero, and hence can
also be eliminated.
The effect of the complete zero-suppressed node reduction can be seen in
\figref{fig:zdd_motiv_3}.
Note, that the construction guarantees canonicity of resulting
diagrams, see~\cite{Ghasemzadeh05}.

\begin{figure}[htb!]
\begin{center}
\subfigure[$a,b,c$]{%
\includegraphics[scale=0.25]{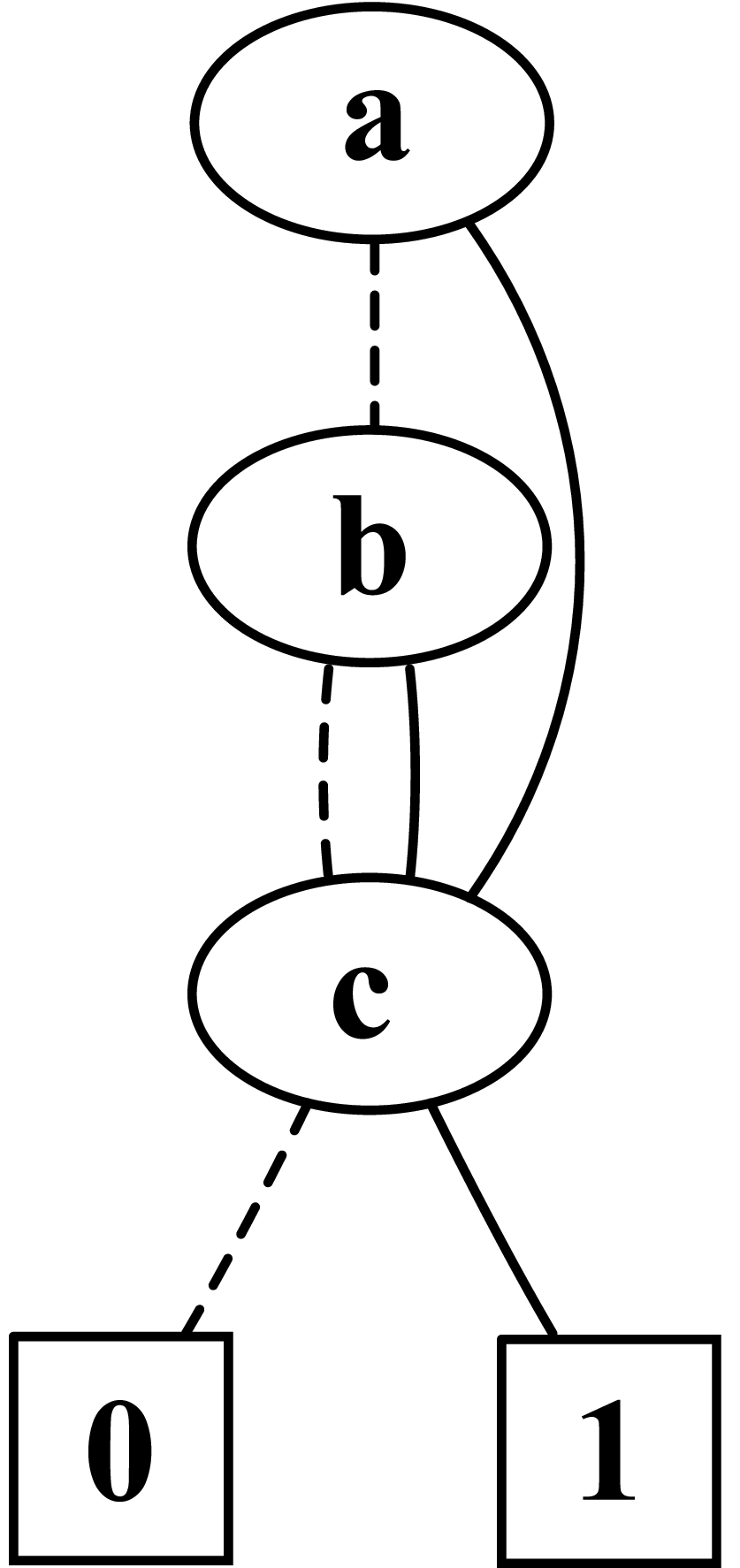}\label{fig:zdd_abc}
}%
\subfigure[$a,c,b$]{%
\hspace{1.5em}%
\includegraphics[scale=0.25]{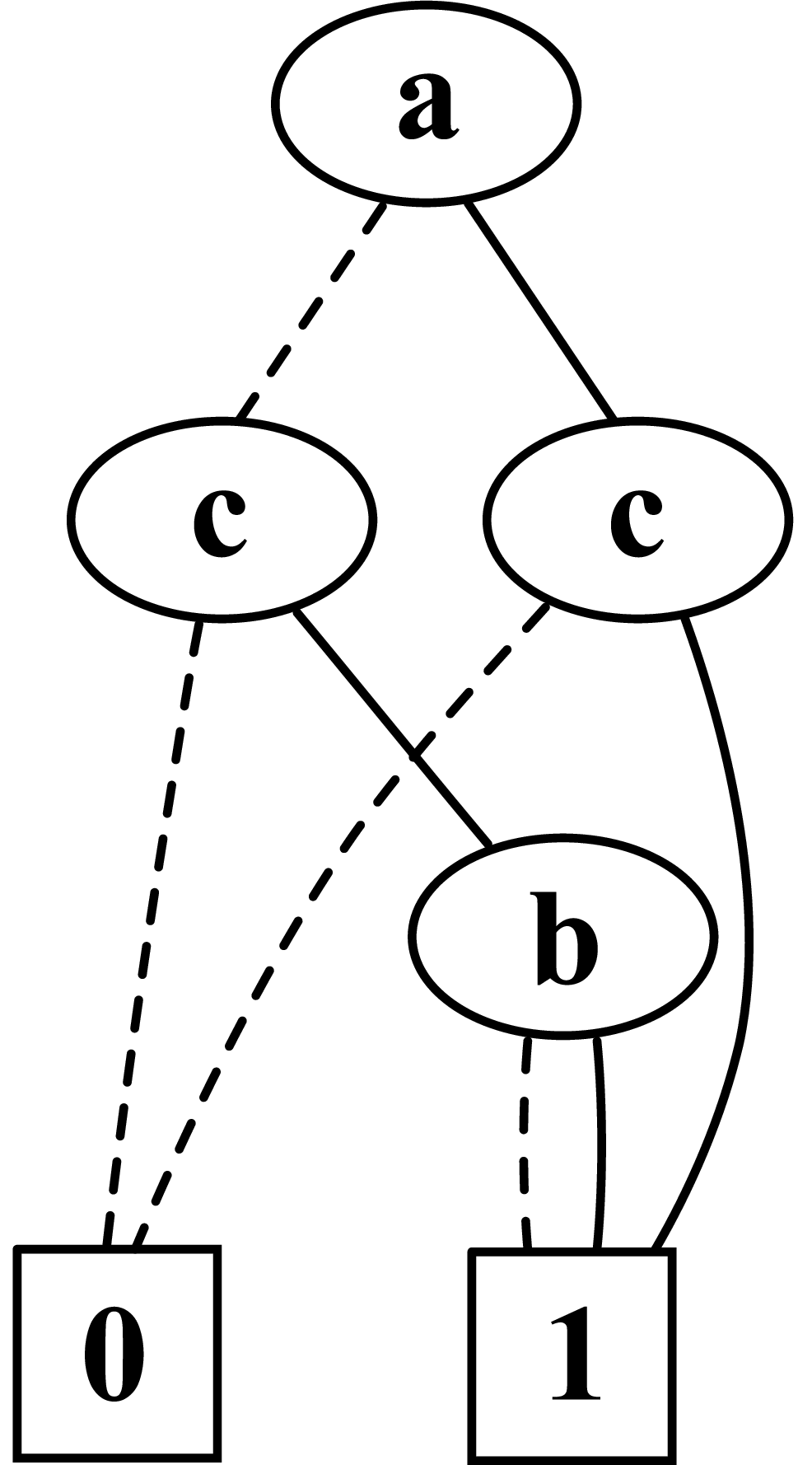}\label{fig:zdd_acb}
}%
\end{center}
\caption{ZDD representing the polynomial~$a\,c+b\,c+c$ for
two different variable orders.
Solid/dashed connections marking then/else-edges, respectively.
\label{fig:zdd_expl}}
\end{figure}

The structure of the resulting ZDD highly depends on the order of the variables,
as \figref{fig:zdd_expl} illustrates.
Hence, a suitable choice of the variable order is always a crucial point,
when \modell{}ing a problem using sets of Boolean polynomials.

Reinterpreting valid paths of a ZDD as terms of a polynomial,
the latter can be accessed in a lexicographical manner, by using the natural
succession arising from the next definition.

\begin{definition}
\label{def:naturalpathseq}
Let~$b$ be a ZDD.
\begin{itemize}
\item
Let~$n_1, n_2,\ldots, n_{m+1}$ be a series of connected nodes
starting at the root node of~$b$ with~$n_{m+1}=1$.
Then the sequence~$(n_1, n_2,\ldots, n_m)$
is called a~\defemph{path} of~$b$.

\item
Let~$x_1>x_2>\ldots>x_n$ be the fixed order of the decision variables.
For two paths~$P=(n_1, n_2,\ldots, n_p)$ and~$Q=(\tilde n_1, \tilde n_2,\ldots,
\tilde n_q)$, the \defemph{natural path ordering}~$<$  is given as:

$P<Q \iff\;\text{there exists a}\; j \in\{1,\ldots, m +1\}, m=\min(p,q)$
such that
\\%
\(
x(n_i) = x(\tilde n_i)\;\text{for}\; 1\leq i<j
\;\text{and}\;
\left\{
\begin{array}{cl}
x(n_j) < x(\tilde n_j) &\text{if}\; j\leq m\\
p < q  &\text{if}\;j=m+1\mcomma
\end{array}
\right.
\)
\\%
where~$x(n)$ denotes the decision variable of a node~$n$.

\item
The ordered sequence~$(P_1, P_2,\ldots, P_s)$
of all paths in~$b$,
is called the
\defemph{natural path sequence} of~$b$.

\end{itemize}
\end{definition}

Note,  that the natural path sequence~$(())$ of the 1-terminal
consists of the empty path only, while path sequence~$()$ of the
0-terminal is empty itself.

One can easily iterate over all paths of a given ZDD. The first path starts at
the root node and follows the~$\thenBranch$ edges, until the 1-terminal is
reached.
For a given path~$P=(n_1, \ldots, n_m)$
the next path in  the natural
path sequence, the \defemph{successor}~$\successor P$ of~$P$, can be computed
follows: let $n_t$ be the first element of~$P$,
with~$\elseBranch(n_i) = 0$, for all~$i > t$, and let
the sequence~$(\tilde n_1, \ldots \tilde n_r)$ denote the first path
in~$\elseBranch(n_t)$, then~$\successor P= (n_1, \ldots, n_{t-1}, \tilde n_1, \ldots \tilde n_r)$.

Although graph-based approaches using decision diagrams for polynomials were
already proposed before, they were not capable of handling
algebraic problems efficiently.
This was mainly due to the fact that the attempts were applied
to very general polynomials, which cannot be represented efficiently as binary
decision diagrams.
For instance, a proposal for utilizing ZDDs for polynomials
with integer coefficients can be found in \cite{Minato95}.
But
Boolean polynomials can be mapped to ZDDs very naturally, since the polynomial
variables
are in one-to-one correspondence with the decision variables in the diagram.
By abuse of notation, we may write in the following~$p$ for the ZDD of a Boolean
polynomial~$p$.

Also, the importance of nontrivial monomial orderings
prevented the use of ZDDs   so far.
In order to enable fast access to leading terms and efficient iterations over
all polynomial terms,
these are usually stored as sorted lists, with respect to a given
monomial ordering \citep{BS:98}.
In contrast, the natural path sequence in  binary decision diagrams is given in
a lexicographical way. Fortunately, it is possible to implement a search for
the leading term and term iterators with moderate effort. Moreover, the results of
basic operations like polynomial arithmetic do not depend on the ordering.
Hence,  these can efficiently be done by using basic set operations.



\subsection{Boolean Polynomial Arithmetic}
\label{sec:poly_arith}

Polynomial addition and multiplication are an essential prerequisite for the
application of \Groebner-based algorithms and related procedures.
In the case of Boolean polynomials, these operations can be implemented as set
operations.
As mentioned in \secref{sec:poly_as_set}, Boolean
polynomials~$p,q\in\Bool$ can be identified with sets~$S_p, S_q
\in \powerset{\powerset{x_1,\ldots,x_n}}$,
\sth~$p =  \sum_{s\in S_p} \left( \prod_{x_\nu \in s}x_\nu \right)$
and~$q =  \sum_{s\in S_q} \left( \prod_{x_\nu \in s}x_\nu \right)$.

Addition  is then just given as~%
$p + q =  \sum_{s\in S_{p+q}} \left( \prod_{x_\nu \in s}x_\nu \right)$,
where~$S_{p+q}$ is computed
as~$S_{p+q} = (S_p\cup S_q) \backslash  (S_p\cap S_q)$.
All three operations~-- union,
complement, and intersection~-- are already available as basic ZDD operations.
For practical applications it is appropriate to
avoid large intermediate sets like~$S_p\cup S_q$ and repeated iterations over
the arguments.
Hence, it is more preferable to have a \specialised addition procedure.
\Algref{alg:recursive_add} below shows a recursive approach for such an
addition.
\begin{algorithm}[htb]
\caption{Recursive addition $h = f+g$}
\label{alg:recursive_add}
\begin{algorithmic}
\REQUIRE $f, g \in \Bool$

\IF{$f = 0$}
\STATE  $h = g$
\ELSIF{$g = 0$}
\STATE  $h = f$
\ELSIF{$f = g$}
\STATE  $h = 0$
\ELSE
\IF{$\mathrm{isCached}(+, f, g)$}
\STATE $h =\mathrm{cache}(+, f, g)$
\ELSE

\STATE set $x_\nu = \mathrm{top}(f)$,  $x_\mu = \mathrm{top}(g)$

\IF{$\nu < \mu$}
\STATE $h = \ite(x_\nu, \thenBranch(f),  \elseBranch(f) + g) $
\ELSIF{$\nu > \mu$}
\STATE $h = \ite(x_\mu, \thenBranch(g), f + \elseBranch(g)) $
\ELSE
\STATE $h = \ite(x_\nu, \thenBranch(f) + \thenBranch(g),
                                \elseBranch(f) + \elseBranch(g) ) $
\ENDIF
\STATE $\mathrm{cache}(+, f, g) = h$
\ENDIF
\ENDIF
\RETURN $h$
\end{algorithmic}
\end{algorithm}

Right after the initial if-statements, which handle trivial cases,
the procedure also includes a cache lookup.
The lookup can be implemented cheaply, because polynomials
have a unique representation as ZDDs.
Hence, previous computations  of the sums of the form~$f+g$ can be reused.
The advantage of a recursive formulation is, that
this also applies to those subpolynomials, which are
generated by~$\thenBranch(f)$ and~$\elseBranch(f)$.
It is very likely, that common subexpressions can be reused during
\Groebner base computation, because of the recurring
multiplication and addition operations, which are used
in Buchberger-based algorithms for elimination of leading terms and the
tail-reduction process.

In a similar manner Boolean multiplication is given in
\algref{alg:recursive_mult}. Note that the procedure computes the unique
representative of the Boolean product~(modulo the field equations). This
multiplication is denoted by~\boolemult\ in the following, while~$\cdot$ means
the usual multiplication. If variables of right- and left-hand side polynomials
are distinct, both operations coincide.
\begin{algorithm}[htb]
\caption{Recursive multiplication $h = f\boolemult g$}
\label{alg:recursive_mult}
\begin{algorithmic}
\REQUIRE $f, g \in \Bool$

\IF{$f = 1$}
\STATE  $h = g$
\ELSIF{$f = 0$ or $g = 0$}
\STATE  $h = 0$
\ELSIF{$g = 1$ or $f = g$}
\STATE  $h = f$
\ELSE
\IF{$\mathrm{isCached}(\boolemult, f, g)$}
\STATE $h =\mathrm{cache}(\boolemult, f, g)$
\ELSE
\STATE  $x_\nu = \mathrm{top}(f)$,  $x_\mu = \mathrm{top}(g)$

\IF{$\nu < \mu$}
\STATE set $p_1 =\thenBranch(f)$, $p_0 =\elseBranch(f)$,
$q_1 = g$, $q_0 = 0$
\ELSIF{$\nu > \mu$}
\STATE set $p_1 =\thenBranch(g)$, $p_0 =\elseBranch(g)$,
$q_1 = f$, $q_0 = 0$
\ELSE
\STATE set $p_1 =\thenBranch(f)$, $p_0 =\elseBranch(f)$,
$q_1 =\thenBranch(g)$, $q_0 =\elseBranch(g)$
\ENDIF
\STATE $h = \ite(x_{\min(\nu, \mu)}, p_0\boolemult q_1 + p_1\boolemult q_1 + p_1 \boolemult q_0,  p_0
\boolemult q_0 ) $
\STATE $\mathrm{cache}(\boolemult, f, g) = h$
\ENDIF
\ENDIF
\RETURN $h$
\end{algorithmic}
\end{algorithm}

\subsection{Monomial Orderings}
\label{sec:monom_order}

While the operations treated in \secref{sec:poly_arith} are independent of the
actual monomial ordering, many operations used in \Groebner
algorithms require such an ordering.
Using ZDDs as basic data structure already yields a natural ordering on Boolean
polynomials as the following theorem shows.

\begin{thm}
Let~$f$ be a Boolean polynomial and~$z$ the corresponding~$ZDD$.
If~$P$ is a path in~$z$, then~$m=\prod_{n_\nu\in P}x(n_\nu)$,
with~$x(n)$ denoting the decision variable of a node~$n$, is a term (and
monomial) in~$f$.
Furthermore, the natural path sequence~$(P_1, P_2,\ldots, P_s)$ yields the
monomials of~$f$ in lexicographical order, and the first path of~$z$ determines
the lexicographical leading monomial of~$f$.
\end{thm}

\begin{pf}
First note, that for a given path~$(n_1,n_2,\ldots, n_m)$,
its ordered sequence of decision variables~$(x(n_1),
x(n_2),\ldots, x(n_m))$ denotes a formal word
in~$x_1,\cdots,x_n$, which can be identified with the monomial given by the
product~$x(n_1)\cdot x(n_2)\cdot\ldots \cdot x(n_m)$.
The first statement is then a consequence of the representation of
polynomials as decision diagrams and the node elimination rule of ZDDs.
The natural ordering of \defref{def:naturalpathseq}
defines then an ordering on the corresponding formal words.
The latter  coincides with the lexicographical ordering,
by comparison of the definitions. Therefore, the natural path sequence
yields the monomials of a polynomials lexicographically ordered, starting with the
leading term.
\end{pf}

Monomials can be represented as single-path ZDDs.
This enables procedures of monomials, analogously to an implementation using
linked lists, but due to the canonicity of the binary decision diagram, equality
check is immediate.
From the implementation point of view, it is not always necessary to generate a
ZDD-based representation for a monomial.
In case, that just some properties are to be checked, and the monomial is not
used in the further procedure, these tests can
also be done on a stacked sequence of nodes, representing a
path in the ZDD.
This kind of stack is used in procedures, which iterate over all terms \wrt
the natural path sequence of a ZDD. Hence, in this case it is already available
without additional costs.

\subsubsection{Degree and block orderings}
\label{sec:deg_ord}

Support of degree orderings are important for \Groebner
algorithms, for two reasons.
First of all, they are necessary for certain algorithms, and second,
because of their better performance in most cases.
%
%
%
A \naive approach would be unrolling all possible paths first, generating all
monomials, and selecting the first among those of maximal degree.
But this procedure could not be cached efficiently.
For a Boolean polynomial~$p = x\cdot p_1 + p_0$ with top
variable~$x$ a recursive formula is 
\begin{equation}
\lm(p) =\left\{
\begin{array}{rl}
 x\cdot \lm(p_1) &\mbox{if}\; \deg(\lm(p_1))+1 \geq\deg(\lm(p_0)) \\
\lm(p_0)& \mbox{else}\mfstop
\end{array}
\right.
\end{equation}
But still this variant accumulates many
single-serving terms.
This can be avoided
by calculating~$\deg(f)=\max(\,\deg(\thenBranch(f))+1, \deg(\elseBranch(f))\,)$
separately.
Caching $\deg(f)$  makes the degree available 
for all recursively generated subpolynomials.
\Algref{alg:recursive_lead_dlex} \utilises this
for computing~$\lm(f)$.
\begin{algorithm}[htb]
\caption{Degree-lexicographical leading term $\lm(f)$}
\label{alg:recursive_lead_dlex}
\begin{algorithmic}
\REQUIRE $f\in\Bool$

\STATE \textbf{if} $\deg(f) = 0$ \textbf{then return} 1
\IF{not $\mathrm{isCached}(\lm, f)$}
\IF{$\deg(f) = \deg(\thenBranch(f)) + 1$}
\STATE $\mathrm{cache}(\lm, f)  = \mathrm{top}(f)\cdot\lm(\thenBranch(f))$
\ELSE
\STATE $\mathrm{cache}(\lm, f)  = \lm(\elseBranch(f))$
\ENDIF
\ENDIF
\RETURN $\mathrm{cache}(\lm, f)$
\end{algorithmic}
\end{algorithm}
Similarly, monomial comparisons and path sequences which yield
polynomial terms in  degree-lexicographical order can be implemented.

A degree-reverse-le\-xi\-co\-gra\-phi\-cal ordering can be handled
in a similar manner. But for this purpose, it is more efficient to reverse the
order of the variables, and the search direction as well.
%
In particular, the leading monomial 
corresponds to \emph{last} path in the natural path
sequence with maximal cardinality, and
\algref{alg:recursive_lead_dlex} can easily
be adapted to this case
by replacing the condition~$(\deg(f) = \deg(\thenBranch(f)) + 1)$
by~$(\deg(f) \not= \deg(\elseBranch(f)))$.

Another important feature are block orderings made of degree orderings.
For this purpose, a block degree can be computed by equipping the
degree-computation with a second argument, which marks the end of the current block~(\ie that block containing the top variable).
Having such a~$\blockdeg$ functionality at hand the leading term computation for
a composition of degree-lexicographical orderings can be obtained by extending
\algref{alg:recursive_lead_dlex} with an iteration over all blocks.


\subsection{Theory of Boolean \Groebner Bases}

In this section,
we present the theory of \Groebner bases over Boolean rings. In the following, we always assume,
that the monomial ordering is global (so $\lm(x^2+x)=x^2$ for every variable $x$).
Since~$\Bool
\cong \Ztwoxoneton/\fieldideal$ this is mathematically equivalent to the theory
of \Groebner bases over the quotient ring. In the classical setting this would
mean to add the field polynomials~\FP to
the given generators~$S \subseteq\Bool$ of a polynomial ideal  and compute a \Groebner basis
of~$\ideal{S, \FP}$ in~$\Ztwoxoneton$.
This general approach is not well-suited for the special case of ideals
representing Boolean reasoning systems.
Therefore,  we propose and develop algorithmic enhancements
and improvements of the underlying theory of \Groebner bases  for
ideals over~$\Ztwoxoneton$ containing the field equations.
Using Boolean multiplication this is implementable directly
via computations with canonical representatives in the quotient ring.
%
%
%
The following theorems shows, that it suffices to treat the Boolean
polynomials introduced in \secref{sec:def_boole_poly} only.

\begin{theorem}
\label{thm:BooleGB}
 Let $S\subseteq \Ztwoxoneton$ be a generating system of some ideal,
such
that~$\fieldequationset \subseteq S \subseteq \Bool \cup \fieldequationset$.
    Then all polynomials created in the classical Buchberger algorithm applied to~$S$ are either Boolean polynomials or field polynomials, if a reduced normal form is used.
\end{theorem}
\begin{pf}
    All input polynomials fulfill the claim.
    Furthermore, every reduced normal form of an s-polynomial is reduced against $\fieldequationset$, so it is Boolean.
    Moreover, using Boolean multiplication every polynomial inside the normal form algorithm is Boolean. Using Boolean multiplication at this point is equivalent to usual multiplication and a normal form computation against the ideal of field equations afterwards.
\end{pf}
\begin{rem}
    Using this theorem we need field equations only in the
generating system and the pair set. On the other hand, we can implicitly
assume, that all field equations are in our polynomial set, and then replace the
pair~$(x_i,p)$ (using Boolean multiplication) by the Boolean polynomial given
as~$x_i\boolemult p=NF(\spoly(x_i,p)|\fieldequationset)$. In this way we can
eliminate the field equations completely. A more efficient implementation would
be to represent the pair by the tuple $(i,p)$, as this still allows the
application of the criteria, but delays the multiplication.
\end{rem}

\begin{lemma}
    The set of field equations $\fieldequationset$ is a \Groebner basis.
\end{lemma}
\begin{pf}
    Every pair of field equations has a standard representation by the product criterion.
    Hence $\fieldequationset$ is a \Groebner basis by Buchberger's Criterion~\citep[Theorem 1.7.3]{singintro}
\end{pf}


\begin{theorem}
    Every $I\subseteq \Ztwoxoneton$ with $I\supseteq \fieldideal$ is radical.
\end{theorem}

\begin{pf}
    Consider $p\in \Ztwoxoneton$, w.\,l.\,o.\,g.\ assume $p$ is reduced against the
leading ideal~$\li(I)$. In particular $\lm(p)$ is a Boolean polynomial.
    Let~$n>0$ and~$q$ be the unique reduced normal form of $p^n$ \wrt the field ideal. So $q$ is also a Boolean polynomial.
    Since $p^n-q$ is a linear combination of field equations,~$p^n-q$ is the zero
function over~$\Ztwo$. By \corref{existszeroone} we get~$p=q$, since~$p^n$
and~$p$ define the same Boolean function.
    Suppose now~$p^n\in I$. Then we have~$p=q=p^n-(p^n-q)\in I$, since~$I\supset \fieldideal$.
\end{pf}

Note that for~$\FP \subseteq I \subseteq \Ztwoxoneton$
the algebraic set~$\V(I)$
is equal to the a priori larger set~$\{\vars \in \overline{\Ztwo\!}\,^n | f(x) = 0 \,\forall f \in I \}$,
 where~$\overline{\Ztwo\!}\,$ denotes the algebraic closure of~$\Ztwo$. Hence
we have

\begin{corollary}
    \label{stronger-nullstellen}
    For ideals $I\subseteq \Ztwoxoneton$ with $I\supseteq \fieldideal$ the
following stronger version of Hilbert's Nullstellensatz holds:
\begin{enumerate}
\item    \( I = \ideal{1} \iff \V(I)=\emptyset\mcomma  \)
\item    \(\I(\V(I))=I\mfstop\)
\end{enumerate}
\end{corollary}


\begin{lemma}
    \label{principal-ideal-lemma}
    If~$I=\pwithfe{p}$ then~$V(I)=V(p)$ and every polynomial $q \in \Ztwoxoneton$ with $\V(q)\supset \V(p)$ lies in $I$.
\end{lemma}
\begin{pf}
    Simple application of Hilbert's Nullstellensatz.
\end{pf}

It is an elementary fact, that systems of logical expressions can be described by a single
expression, which describes the whole system \behaviour.
Hence, the one-to-one correspondence of Boolean polynomials and Boolean
functions given by the mapping defined in \defref{def:maplogic} motivates the
following theorem.

\begin{theorem}
    \label{unique-generator}
    Every ideal in $\Ztwoxoneton/\fieldideal$ is generated by the equivalence class of one unique Boolean polynomial.
In particular,~$\Ztwoxoneton/\fieldideal$ is a principal ideal ring~(but not a
domain).
\end{theorem}

\begin{pf}
    We use the one-to-one correspondence of ideals in the quotient ring and ideals in $\Ztwoxoneton$ containing $\fieldideal$.
    Therefore, let $\fieldideal \subset I\subset \Ztwoxoneton$.
    By \corref{cor:boole1-1} there exists a Boolean polynomial~$p$
\sth~$\V(\pwithfe{p})=\V(I)$.
    By \thmref{stronger-nullstellen} we get~$I=I(\V(\pwithfe{p}))=\pwithfe{p}$.
    %
    Suppose, there exists a second Boolean polynomial~$q$
with~$I=\pwithfe{q}$. Then $$V(p)=V(I)=V(q).$$
    So $p$ and $q$ define the same characteristic function, which means that they are identical  Boolean polynomials.

\end{pf}

Hence, using \thmref{polys-function},    \corref{cor:boole1-1} and \corref{stronger-nullstellen}, we have the following bijections:
\begin{equation*}
\begin{array}{c}
 \Bool
 \leftrightarrow
 \{\text{Boolean functions}\}
 \leftrightarrow \\
  \{\text{ideals}\; I \subseteq\Ztwoxoneton\;\text{with}\;\FP \subseteq I \}
  \leftrightarrow \\
 \{\text{algebraic subsets of}\; \Ztwo^n\}
\leftrightarrow
 \{\text{subsets of}\; \Ztwo^n\} \mfstop
\end{array}
\end{equation*}

\begin{defn}

       For any subset~$H\subseteq \Ztwoxoneton$, call
       $$BI(H):=\pwithfe{H}\subseteq \Ztwoxoneton$$ the {\defemph{Boolean ideal of H}}.
     We call a reduced \Groebner basis of $\BI(H)$ the {\defemph{Boolean
    \Groebner basis}} of~$H$, short $\BGB(H)$.

\end{defn}

Recall from \thmref{thm:BooleGB} that~$\BGB(H)$ consists of Boolean polynomials and can be extended to a reduced \Groebner basis of $\BI(H)$ by
adding some field polynomials.

\begin{theorem}
    Let $p, q\in\Bool$ with $\V(p)\subset \V(q)$. Then $\pwithfe{p}
\supset\pwithfe{q}$ and we say~$p$ implies~$q$. This implication relation forms a partial order on the set of Boolean polynomials.
\end{theorem}
\begin{pf}
    Since both ideals are radical, Hilbert's Nullstellensatz gives the ideal containment.
    The implication is a partial order by the one-to-one correspondence between
Boolean polynomials and sets. It corresponds itself to the inclusion of sets.
\end{pf}

\subsection{Criteria}

Criteria for keeping the set of critical pairs in the Buchberger algorithm small
are a central part of any \Groebner basis algorithm aiming at practical
efficiency. In most implementations the chain criterion and the product criterion  or variants of them are used.

These criteria are of quite general type, and it is a natural question, whether we can
formulate new criteria for Boolean \Groebner bases. Indeed,  this is the case. There are two types of pairs to
consider: Boolean polynomials with field equations, and pairs of Boolean polynomials.
We concentrate on the first kind of pairs here.

\begin{thm}
	Let $f\in\Bool$ be of the form~$f=l\cdot g$, $l$ a polynomial with
linear leading term~$x_i$, and~$g\in\Ztwoxoneton$ be any polynomial. Then $\spoly(f,x_i^2+x_i)$ has a nontrivial $t$-representation against the system consisting of~$f$ and the field equations.
	\label{linear-factor-criterion}
%
\end{thm}

The theorem was proved by Brickenstein
in~\cite{BrickensteinDreyerPreMega07}.

\begin{lemma}
    \label{trivial-multiplication}
    Let $G$ be a \Groebner basis, $f$ a polynomial, then $\set{f\cdot g\vert g \in G}$ is \Groebner basis.
\end{lemma}

\begin{rem}
    This lemma is trivial, we just want to show the difference to the next theorem.
\end{rem}

\begin{thm}
    Let $G$ be a Boolean \Groebner basis, $l\in\Bool$ with~$\deg(\lm(l))=1$ and
$\support(l)\cap\support(g)=\emptyset$ for all~$g\in G$. Then $\set{l\cdot g\vert g\in G}$ is a Boolean \Groebner basis~%
that is,  $\set{l\cdot g\vert g\in G} \cup \fieldequationset$ is a \Groebner
basis. In other words, we get a \Groebner basis again by multiplying the Boolean polynomials, but not the field equations with the special polynomial $l$.
\end{thm}

\begin{pf}
    We show, that every s-polynomial has a non-trivial $t$-representation.
    We have to consider three types of pairs.
    If~$p$, $q$ are  both field polynomials, $\spoly(p,q)$ has a standard
representation by the product criterion.
    If $p$, $q$ are both Boolean polynomials, then~$\spoly(l\cdot p,l\cdot q)$ has a standard representation by multiplying the standard representation of $\spoly(p,q)$ by $l$.
    Now let $p$ be a Boolean polynomial and $q$ a field polynomial, say $q=x^2+x$.
    If~$\lm(l)=x$, then $\spoly(l\cdot p, q)$ has a nontrivial~$t$-representation by \thmref{linear-factor-criterion}.
    If $x$ occurs in $\lm(p)$, then by \lemref{trivial-multiplication}
$\spoly(l\cdot p,\l \cdot q)$ has a standard representation against~$\set{l \cdot g|g \in G} \cup \set{l\cdot e\vert e \in \fieldequationset}$, so also against
     the set~$\set{l \cdot g\vert g \in G} \cup \fieldequationset$.
    Hence, we just have to show, that the difference
to~$\spoly(l\cdot p,\l \cdot q)$
has a~$t$-representation with~$t<\lm(p)\cdot\lm(l)\cdot x:=c.$
Setting    \[h:=\spoly(l\cdot p,l\cdot (x^2+x))-\spoly(l\cdot p,x^2+x)=\tail(l) \cdot (x^2+x)\]
    we get that~$x^2+x$ divides $h$,
and~$\lm(h)= \lm( (x+1) \cdot \tail(l)) \cdot x < c$,
 since~$\lm(p)$ contains~$x$.
So $h$ has standard representation against~$x^2+x$.
    If $x$ does  neither occur  in $\lm(f)$ nor in $\lm(l)$ the product criterion applies.
    Reducedness follows from the fact, that $l$ does not share any variables with $G$.
\end{pf}

\subsection{Symmetry and Boolean \Groebner bases}

In this section we will show how to use the theory presented in the previous
section to build faster algorithms by using symmetry and simplification by pulling out factors with linear leads.

    For a polynomial $p$ we denote by $\varsof(p)$ the set of variables actually occurring in the polynomial.

\begin{definition}
    Let $p$ be a polynomial in $\Ztwoxoneton$ with a given monomial ordering~$>$,
    $\vert \varsof(p)\vert=k$, $I=\varsof(p)=\set{x_{i_1},\ldots, x_{i_k}}$,
and~$J=\set{x_{j_1},\ldots, x_{j_k}}$ be any set of~$k$ variables.
    We call a morphism of polynomials algebras over~$\Ztwo$,
\[f:\Ztwo [I]\rightarrow
\Ztwo[J]:x_{i_s} \mapsto x_{j_s} \;\text{for all}\; s\mcomma\]
 a suitable shift for $p$, if and
only if for all monomials~$t_1,t_2\in \Ztwo[I]$ the
relation~$t_1>t_2\Longleftrightarrow f(t_1)>f(t_2)$ holds.
\end{definition}

\begin{remark}
    In the following we concentrate on the problem of calculating $\BGB(p)$ for
one Boolean polynomial $p$ (non-trivial, as field equations are implicitely included). So, if
we know $\BGB(q)$ for a Boolean polynomial $q$ and if there exists a suitable
shift $f$ with $f(q)=p$, then $f(\BGB(q))=\BGB(p)$.
Hence, we can avoid the computation of~$\BGB(p)$.
    Adding all elements of $BGB(p)$ to our system means that we can omit all pairs of the form $(p,x_i^2+x_i)$.
    A special treatment (using caching and tables) of this kind of pairs is a
good idea, because this is a often reoccurring phenomenon. As these pairs depend
only on $p$ (the field equations are always the same), this  reduces the number
of combinations significantly.
\end{remark}
\begin{remark}
    Note, that the concept of Boolean Gr\"obner bases fits very well here, as $\BGB(p)$
    is the same in $\Ztwo[\varsof(p)]$ as in $\Ztwoxoneton$, although the last case refers to a \Groebner basis with more field equations.
\end{remark}

\begin{definition}
    We define the relation $p \sim_{pre} q$, if and only if there exists a
suitable shift between $p$ and $q$ or if there  exists an~$l$ with $\deg(\lm(l))=1$ and $p=l\cdot q$.
    From $\sim_{pre}$ we derive the relation $\sim_{sym}$ as its reflexive, symmetric, transitive closure (the smallest equivalence relation containing $\sim_{pre}$).
\end{definition}
\begin{remark}
    For all $p$ and $q$ in an equivalence class of $\sim_{sym}$ the Boolean
\Groebner basis~$\BGB(p)$ can be mapped to $\BGB(q)$ by a suitable variable shift and pulling out (or multiplying) by Boolean polynomials with linear lead.
    In practise, we can avoid complete factorizations by restricting ourselves to detect factors of the form $x$ or $x+1$.
%
Using these techniques it is possible to avoid the explicit calculation of many critical pairs.
\end{remark}

\begin{definition}
    A monomial ordering is called symmetric, if the following holds.
For every~$k$, and every two subsets of variables~$I=\set{x_{i_1},\ldots,x_{i_k}}$,
and~$J=\set{x_{j_1},\ldots,x_{j_k}}$ with $i_z<i_{z+1}$, $j_z<j_{z+1}$
for all~$z$ the~$\Ztwo$-algebra homomorphism
    \[f:\Ztwo[I]\rightarrow \Ztwo[J]:x_{i_z}\mapsto x_{j_z}\] defines a suitable shift.
\end{definition}

\begin{algorithm}
\caption{Calculating $\BGB(p)$ in a symmetric order}
\label{BGB}
\begin{algorithmic}
\REQUIRE $p\in \Bool$, $>$ a monomial ordering
\ENSURE  $\BGB(p)$
\STATE{pull out as many factors with linear lead as possible}
\STATE{calculate a more canonical representative $q$ of the equivalence class of $p$ in $\sim_{sym}$ by shifting $p$ to the first variables}
\IF{$q$ lies in a cache or table}
    \STATE{$B := \BGB(q)$ from cache}
\ELSE
\STATE{$B := \BGB(q)$ by Buchberger's algorithm}
\ENDIF
\STATE{shift $B$ back to the variables of $p$}
\STATE{multiply $B$ by the originally pulled out factors}

\RETURN $B$
\end{algorithmic}
\end{algorithm}

For a symmetric ordering it is always possible to map a polynomial $p$ to the
variables $x_1,\ldots,x_{\vert \varsof(p) \vert}$ by a suitable shift.
This is \utilised in \algref{BGB} for speeding up calculation of Boolean \Groebner
bases.
    In the following we assume that the representative chosen in the algorithm is canonical (in particular uniquely determined in the equivalence class in $\sim_{sym}$), if every factor with linear lead is pulled out.
\begin{remark}
    From the implementation point of view, it turned out to be useful to
 store the $\BGB$ of all $2^{16}$ Boolean polynomials in up to four variables in a precomputed table, for more variables we use a dynamic cache (pulling out factors reduces the number of variables).
Using canonical representatives increases the number of cache hits.

The technique for avoiding explicit calculations can be integrated in nearly every algorithm similar to the Buchberger's algorithm.
Best results were made by combining these techniques with the algorithm slimgb
\citep{SlimBa06}, we call this combination $\symmgbGFTwo$. For our computations the strategy in slimgb for dealing with elimination orderings is quite essential.
\end{remark}

\subsubsection*{Practical meaning of symmetry techniques}
The real importance of symmetry techniques should not only be seen in avoiding
computations in leaving out some pairs. In constrast,
application of the techniques described above 
changes  the \behaviour of the algorithm  completely.
Having a Boolean polynomial $p$, the sugar value \citep{GiovMNRT:91} of the pair $(p,x^2+x)$ is usually $\deg(p)+1$, which corresponds to the position in the waiting queue of critical pairs.
It often occurs that in $\BGB(p)$ polynomials with much smaller degree occur.

Having these polynomials earlier, we can avoid many other pairs in higher degree.
This applies quite frequentely in this area, in particular, when we have many
variables, but the resulting \Groebner basis looks quite simple~(for example
linear polynomials). The earlier we have these low degree polynomials, the
easier the remaining computations are, resulting in less pairs and faster normal form computations.
\section{Applications}

The algorithms described in \secref{sec:rgb} resp.~\ref{sec:bgb} have been implemented
in \Singular~\cite{GPS05} resp.\ the \PolyBoRi framework \cite{BrickensteinDreyerPreMega07}. We use these
implementations to test our approach by computing realistic examples from formal verification.
We compare the computations with other computer algebra system and with SAT-solvers, all considered
to be state-of-the-art in their field.

Moreover, we state open questions and conjectures, in particular in the case of \Groebner bases over rings,
an area which is not very much explored.

The application of \Groebner bases over $\Ztwon$ is still under development. Here we mention
mainly problems in connection with the proposed applications. On the other hand we show
that the improvements developed in \secref{sec:nf} and \secref{sec:comp_std_bases} for
\Groebner bases over weak factorial principal rings are extremely useful for computations
over these rings.

\subsection{Standard bases over rings}
Let us recapitulate the original problem first, which was posed in \secref{sec:wordlevel}.

\begin{problem}
Given a finite set of polynomials $\{f_i\}\subset\Ztwon$. Does a common
zero of the system $\{f_i=0\}$ exist, \ie is $\V(\skalar{f_i}) \neq \emptyset$?
\end{problem}

To answer this question with the help of computer algebra and \Groebner bases theory, the following
key problems have to be solved.

\begin{problem}\label{prob:agenda}\
\begin{enumerate}
\item \label{enum:computeRingGB}
An efficient algorithm\footnote{Here and in the following efficient refers to practical performance
and not to the complexity of the algorithms.} to compute \Groebner bases over $\Ztwon$.
\item A way to handle vanishing polynomials, \ie polynomials evaluating to zero everywhere.\label{enum:vanIdeal}
\item A suitable Nullstellensatz equivalent for $\Ztwon[\vars]$, or at least a simple \Groebner basis
criterion for the existence of a common zero over some extension ring.\label{enum:Nullensatz}
\end{enumerate}
\end{problem}

In \secref{sec:rgb} we explained, how an efficient algorithm for \probref{prob:agenda}\eqref{enum:computeRingGB}
can be instantiated. In order to \optimis{}e the algorithm in the case of $\Ztwon$ we can replace all greatest common divisor computations by fast divisibility tests.

We implemented the algorithm in the kernel of the computer algebra system \Singular \cite{GPS05} and compared the performance to Magma, the only other system we found to be capable of computing \Groebner bases in $\Ztwon$. As we could not solve industrial-sized problems due to time and space explosion we compared the implementations with random instances. In  \tabref{table:bench_ringgb2n} we present only a few concrete runtimes, but they give an overall impression of the data.
The table shows that the special algorithms for $Z_m$ (apparently not contained in Magma) pay off substantially.

{
\providecommand{\secunit}{\ensuremath{\,\mathrm{s}}}
\providecommand{\MBunit}{\ensuremath{\,\mathrm{MB}}}

\providecommand{\timedoutAfter}[1]{\multicolumn{1}{c}{$\infty$} & }
\providecommand{\memoryout}{ &\,\hfill$\infty$\hfill\, }
\providecommand{\timedout}{\multicolumn{2}{c}{time out after 1h}}

\begin{table}
\begin{center}
{\scriptsize
\begin{tabular}{cccc|c|rr|rr}
\#vars.&\#polys. &maxdeg&$\frac{\text{\#mons.}}{\text{\#polys.}}$ & \#GB &
    \multicolumn{2}{c|}{\Singular}  &  \multicolumn{2}{c}{Magma}  \\\hline

2 & 5 & 15 & 69.2 & 3 & 0.40\secunit & 4.11\MBunit & 68.16\secunit & 13.57\MBunit \\
3 & 3 & 10 & 6.7 & 254 & 8.50\secunit & 17.23\MBunit & 1287.80\secunit & 19.60\MBunit \\
3 & 3 & 15 & 7.4 & 599 & 204.82\secunit & 146.98\MBunit & \timedout \\
4 & 4 & 10 & 2.8 & 120 & 0.04\secunit & 0.87\MBunit & 10.68\secunit & 9.52\MBunit \\
4 & 4 & 10 & 3.0 & 361 & 20.36\secunit & 32.24\MBunit & \timedout \\
5 & 5 & 10 & 2.4 & 584 & 0.15\secunit & 1.09\MBunit & 455.35\secunit & 30.07\MBunit \\
5 & 5 & 10 & 2.8 & 1043 & 1.11\secunit & 2.34\MBunit & \timedout \\
7 & 5 & 10 & 2.0 & 614 & 0.14\secunit & 1.14\MBunit & 40.06\secunit & 35.35\MBunit \\
7 & 5 & 10 & 2.2 & 2547 & 2.23\secunit & 3.03\MBunit & \timedout \\
10 & 10 & 4 & 1.9 & 436 & 0.11\secunit & 1.09\MBunit & 92.45\secunit & 16.75\MBunit \\
10 & 10 & 4 & 3.0 & 11734 & 963.39\secunit & 341.70\MBunit & \timedout \\
12 & 10 & 3 & 2.3 & 5536 & 18.40\secunit & 16.75\MBunit & \timedout \\
12 & 10 & 3 & 3.0 & 1940 & 3.69\secunit & 13.12\MBunit & \timedout
\end{tabular}
}
\end{center}
\caption{Computation of a Gr\"{o}bner basis in $\Z_{2^{10}}$ with degree reverse lexicographical ordering.
Randomly generated examples on an AMD Dual Opteron 2.2 GHz, 16 GB RAM.}
\label{table:bench_ringgb2n}
\end{table}
}

To deal with \probref{enum:vanIdeal}, that is with the ideal of vanishing polynomials in~$\Z_{m}$ with~$m\in\N$
we determined the minimal \Groebner basis $G_0$ of
$$I_0 := \{f\in\Z_m~|~\forall\vars:f(\vars)=0\}$$
combinatorially (cf. \cite{vangb-oli}). The size of $G_0$ grows roughly with $\smaranch(m)^{\# \text{variables}}$, where $\smaranch(m)$ is the
Smarandche function \cite{gen_smarand}. Hence, for a typical application instance of formal verification just
listing the ideal $G_0$ becomes infeasible. We therefore devised a method of constructing only the necessary elements
of $G_0$ for $s$-polynomial and normal form computations, but even their number grows exponentially in
the number of variables.

Another obstacle, related to this one, arises while investigating the modeling strength of polynomials functions in comparison to arbitrary functions
from $\Z_m^n\to\Z_m$. Here we have the following

\begin{observation}[\cite{vangb-oli}]
There are many more functions~$\Z_m^n\to\Z_m$ than polynomial functions and many more
subsets of~$\Z_m^n$ than varieties if~$m$ is not a prime number. The quotient of all functions by polynomial functions grows
at least double-exponentially in the number of variables.
If~$m$ is a prime, then all functions respectively subsets of $\Z_m^n$ are polynomial, respectively algebraic.
\end{observation}

The following conjecture was verified for small $m, n$.

\begin{conj}
A function $\Z_m^n\to\Z_m$ is polynomial if and only if Newton interpolation works. This
means that the division during the algorithm is possible, but not necessarily unique.
\end{conj}


With respect to \probref{prob:agenda} \eqref{enum:Nullensatz} we mention
the following lemma which is a negative result.

\begin{lem}\label{lem:no_alg_closure}
Let $\cR$ be a ring with zero divisors. There exists no ring $\hat{\cR}\supset \cR$, such that every
non-constant polynomial of $\cR[x]$ has a zero in $\hat{\cR}$.
\end{lem}

\begin{proof}
Let $n\in\cR\wo\{0\}$ be a zero divisor and consider $f=n x - 1$. Assume there exists a ring $\hat{\cR}\supset\cR$ which contains a root $r$
of $f$.
Then
$f(r)=n\cdot r - 1 = 0$ and hence $1=n\cdot r$.
On the other hand, there exists an $m\neq 0$ with $m\cdot n=0$ and hence $m\cdot 1 = m\cdot n \cdot r=0$, a contradiction.
\end{proof}

\begin{remark}
If $C$ has no zero divisors then a ring $\hat{\cR}$ as in \lemref{lem:no_alg_closure} exists. We may take $\hat{\cR}$ just
as the algebraic closure of the quotient field of $\cR$.
If $I$ is an ideal in $\cR[\vars]$ we set $\hat{\V}(I):=\{\vars\in \hat{\cR}^n~|~f(\vars)=0~\forall f\in I\}$
and get the following answer to \probref{prob:agenda} \eqref{enum:Nullensatz}: Let $G\subset \cR[\vars]$ be
a \Groebner basis of $I$. Then $\hat{\V}(I)=\emptyset$ iff $G$ contains a non-zero element of $\cR$.

However, if $\cR$ has zero divisors, it is not clear how a useful answer to \probref{prob:agenda} \eqref{enum:Nullensatz}
should look like.
\end{remark}


\subsection{The \PolyBoRi Framework}
We will give a brief description of the \PolyBoRi
framework~\citep{BrickensteinDreyerPreMega07} and the implemented
algorithms. At the end of this section, the time and space requirements of some
benchmark examples are compared with
those of other computer algebra systems and a SAT-solver.


The core routines of \PolyBoRi form a C++ library for
\emph{Poly}nomials over \emph{Bo}olean \emph{Ri}ngs
providing high-level
data types for Boolean polynomials and monomials, exponent vectors, as well as
for the underlying Boolean rings.
%
The ZDD structure, which is used as internal storage for
polynomials and monomials,
is based on a data type from 
CUDD~\citep{somenzi05cudd}.

In addition, basic polynomial operations~-- like addition and multiplication~--
have been implemented and associated to the corresponding operators.
\PolyBoRi's polynomials also provide ordering-dependent functionality, like
lead\-ing-term computations, and iterators for accessing polynomial terms in the
style of \emph{Standard Template Library}'s iterators
\citep{stepanov94standard}.
This is implemented by a stack, which holds a valid path.
The corresponding monomial may be returned on user request, and
incrementing the iterator results in a search for a valid path, corresponding to
next term in monomial order.
The ordering-dependent functions are currently
available for the orderings introduced in \secref{sec:monom_order}
and block orderings thereof.

Issues regarding the monomial ordering
and the internal data structure are hidden behind
a user programming interface. This allows the formulation of generic
procedures in terms of computational algebra, without the need for caring
about internals.
This will then work for any applicable and implemented Boolean ring.

Complementary, a complete Python \citep{RossumPython01} interface allows
parsing of complex polynomial systems.
Rapid prototyping of sophisticated and
easy extendable strategies
for \Groebner base computations was possible by using this script language.
%
With the tool ipython the \PolyBoRi data structures and procedures can
be used interactively.
In addition, interfaces to the computer algebra system
\Singular{}~\citep{GPS05} und the \Sage system~\citep{sage} are under development.

\subsection{Timings}
\providecommand{\secunit}{\ensuremath{\,\mathrm{s}}}
\providecommand{\MBunit}{\ensuremath{\,\mathrm{MB}}}

\providecommand{\timedoutAfter}[1]{\multicolumn{1}{c}{$\infty$} & }
\providecommand{\memoryout}{ &\,\hfill$\infty$\hfill\, }
\providecommand{\timedout}{\timedoutAfter{1\,h}}

This section presents some benchmarks comparing \PolyBoRi to general purpose and
\specialised computer algebra systems.
The following timings have been done on a AMD Dual Opteron 2.2 GHz (all systems
have used only one CPU) with 16 GB RAM on Linux.
The used ordering was lexicographical, with the exception of FGb, where
degree-reverse-lexicographic was used.
\PolyBoRi also implements degree
orderings, but for the presented practical examples elimination orderings seem to
be more appropriate. A recent development in \PolyBoRi was the implementation of
block orderings, which behave very natural for many examples.

We compared the computation of a Gr\"{o}bner basis for the following system releases with the development version
of \PolyBoRi's \symmgbGFTwo:

\begin{tabular}{lll}
Maple& 11.01, June 2007 & Gr\"obner package, default options\\
FGb &  1.34, Oct.\ 2006 &  via Maple 11.01, command: fgb\_gbasis\\
\Magma & 2.13-10, Feb.\ 2007&  command: GroebnerBasis, default options\\
\Singular& 3-0-3, May 2007 & std, option(redTail)\\
\end{tabular}

 Note, that this  presents the state of
\PolyBoRi in the development version  in August 2007 only. Since the
project is very young there is still room for major performance
improvements.
The examples were chosen from current research problems in formal verification. %
%
%
All timings of the computations (lexicographical ordering) are \summarised in \tabref{table:results}.
%
\begin{table}
\begin{center}
{\scriptsize
\newcommand{\cstyle}{{\hspace{.6em}}}
\renewcommand{\MBunit}{}
\renewcommand{\secunit}{}
\begin{tabular}{@{}crr|rr||@{\cstyle}r@{\cstyle}r@{\cstyle}|@{\cstyle}r@{\cstyle}r@{\cstyle}|@{\cstyle}r@{\cstyle}r@{\cstyle}|@{\cstyle}r@{\cstyle}r@{\cstyle}}
\multicolumn{3}{c}{} & \multicolumn{2}{c}{\PolyBoRi} &
  \multicolumn{2}{c}{FGb}&
  \multicolumn{2}{c}{Maple}&
  \multicolumn{2}{c}{\Magma} & \multicolumn{2}{c}{\Singular}
\\
Example &\multicolumn{2}{@{\cstyle}c|}{Vars./Eqs.} &

\multicolumn{1}{@{\cstyle}c}{s} & \multicolumn{1}{@{\cstyle}c||@{\cstyle}}{MB} &
\multicolumn{1}{@{\cstyle}c}{s} & \multicolumn{1}{@{\cstyle}c|@{\cstyle}}{MB} &
\multicolumn{1}{@{\cstyle}c}{s} & \multicolumn{1}{@{\cstyle}c|@{\cstyle}}{MB} &
\multicolumn{1}{@{\cstyle}c}{s} & \multicolumn{1}{@{\cstyle}c|@{\cstyle}}{MB}  &
\multicolumn{1}{@{\cstyle}c}{s} & \multicolumn{1}{@{\cstyle}c}{MB}

\\
\hline
mult4x4&  55 &  48  &

 0.00\secunit & 54.54\MBunit
& 
1.76\secunit & 5.50\MBunit
 & 
1.96\secunit & 4.87\MBunit
 & 
0.91\secunit & 10.48\MBunit
 & 
0.02\secunit & 0.66\MBunit
\\
mult5x5&  83 &  74  &

 0.01\secunit & 54.66\MBunit
& 
219.09\secunit & 6.37\MBunit
 & 
236.14\secunit &  6.87\MBunit
 & 
31.28\secunit & 46.05\MBunit
 & 
0.01\secunit & 1.67\MBunit
\\
mult6x6&  117 &  106  &
0.03\secunit & 54.92\MBunit
& 
\multicolumn{2}{c|@{\cstyle}}{failed}%
 & 
\timedout 
 & 
\timedout
 & 
4.28\secunit & 21.19\MBunit
\\
mult8x8&  203 &  188  &

0.40\secunit & 55.43\MBunit
& 
\timedout
 & 
\timedout 
 & 
\timedout
 & 
\memoryout   
\\
mult10x10&  313 &  294  &

 18.11\secunit & 85.91\MBunit
& 
\memoryout
 & 
\memoryout & 
\memoryout & 
\memoryout 
\end{tabular}
}
\end{center}
\caption{Timings and memory usage for benchmark examples. The~$\infty$ symbols
in time and memory columns
mark timeout after~1 hour and out of memory at~15\,GB.}
\label{table:results}
\end{table}
%

The authors of this article are convinced, that the default strategy of
\Magma\ is not well suited for these examples (walk, see \cite{CoKaMa97}, or
\homogenisation). However, when we tried a direct approach in \Magma, it ran
very fast out of memory~(at least in the larger examples). We can conclude,
that the implemented \Groebner basis algorithm in \PolyBoRi offers a good
performance combined with suitable memory consumption.
Part of the strength in directly computing \Groebner bases (without walk or
similar techniques) is inherited from the slimgb algorithm in \Singular.  On the
other hand our data structures provide a fast way to rewrite polynomials, which
might be of bigger importance than sparse strategies in the presented examples.


In order to treat classes of examples, for which the lexicographical ordering is
not the best choice, \PolyBoRi is also equipped with other monomial orderings.
Although its internal data structure is ordered lexicographically,
the computational overhead of degree orderings is small
enough such that the advantage of these orderings come into effect.
\Tabref{tab:orderbench} illustrates this for a series of randomly generated
unsatisfiable uniform
examples~\cite{satlib00}. The latter arise from benchmarking SAT-solvers, which can handle them
very quickly, as their conditions are easy to contradict. But they are still a
challenge for the algebraic approach.
\begin{table}
\begin{center}
{\tiny
\begin{tabular}{cccc|rr||rr|rr}
\\
Example &\multicolumn{2}{c}{Vars./Eqs.} & Order. &

  \multicolumn{2}{c||}{\PolyBoRi} &
  \multicolumn{2}{c|}{\Magma} &
  \multicolumn{2}{c}{FGb}

\\
\hline
uuf50\_10&  50 &  218  &
lp &
8.76\secunit & 71.98\MBunit
 & 
9.77\secunit & 28.21\MBunit
   \\
&&& dlex&
8.98\secunit & 72.53\MBunit
 & 
10.35\secunit & 32.71\MBunit
   \\
 &&& dp\_asc&
8.14\secunit & 72.24\MBunit
 & 
8.40\secunit & 27.42\MBunit
 & 
74.76\secunit & 6.75\MBunit
\\
uuf75\_8&  75 &  325  &
lp &
843.38\secunit & 819.80\MBunit
 & 
14015.21\secunit & 1633.62\MBunit
   \\
&&& dlex &
553.43\secunit & 490.86\MBunit
 & 
14291.45\secunit & 2439.53\MBunit
   \\
 &&& dp\_asc&
448.53\secunit & 472.04\MBunit
 &  
13679.42\secunit & 2539.24\MBunit
& 
\!\!\!99721.46\secunit    & \!\!\!\!\!\!8958.36\MBunit
\\
uuf100\_01&  100 &  430  &
lp &
44779.77\secunit & 12309.79\MBunit
 & 
 \timedoutAfter{2\,days}
\\
  &&&dlex &
11961.86\secunit & 6101.43\MBunit
& 
 \timedoutAfter{2\,days}
\\
&&&dp\_asc &
10635.72\secunit & 6146.47\MBunit
 & 
\timedoutAfter{2\,days}
& 
\multicolumn{2}{c}{failed}
\end{tabular}
}
\end{center}
\caption{Timings and memory usage for Gr\"{o}bner basis computations w.\,r.\,t.\ various
orderings.
The~$\infty$ symbols means timeout after~2 days, \emph{failed} stopped with
error message, and \emph{dp\_asc} denotes \emph{dp} with reversed variable
order.
}
\label{tab:orderbench}
\end{table}
The strength of \PolyBoRi is visible in the more complex examples, as it scales
better than the other systems in tests.

In addition the performance of \PolyBoRi is compared  with the freely
available SAT-solver MiniSat2~(release date 2007-07-21),
which is
state-of-the-art
among publicly available solvers~\citep{een03minisat}.
The examples consist of formal verification
examples corresponding to digital circuits with $n$-bitted multipliers
and the pigeon hole benchmark, which is a
standard benchmark problem for
SAT-solvers, \eg used in in~\cite{satlib00}%
.
The latter checks whether it is possible to place~$n+1$ pigeons in~$n$ holes
without two of them being in the same hole~(obviously, it is unsatisfiable).

%
\begin{table}
\begin{center}
{\scriptsize
\begin{tabular}{l|rr|rr|rr}
& \multicolumn{2}{c|}{Vars./Eqs.}  &
    \multicolumn{2}{c|}{\PolyBoRi}  &  \multicolumn{2}{c}{MiniSat}  \\
\hline
 hole8 & 72 & 297 &
 1.88\secunit &56.59\MBunit &
 0.30\secunit & 2.08\MBunit \\
 hole9 & 90 & 415 &
 8.01\secunit & 84.04\MBunit &
 2.31\secunit & 2.35\MBunit \\
 hole10 & 110 & 561 &
 44.40\secunit & 97.68\MBunit &
 25.20\secunit & 3.24\MBunit \\
 hole11 & 132 &738  &
 643.14\secunit & 130.83\MBunit &
 782.65\secunit & 7.19\MBunit  \\
hole12 & 156 &949   &
  10264.92\secunit & 338.66\MBunit &
   22920.20\secunit & 17.13\MBunit
\\
\hline
 mult4x4 &55 & 48 &
 0.00\secunit & 54.54\MBunit &
 0.00\secunit & 1.95\MBunit \\
 mult5x5 & 83 & 74 &
 0.01\secunit & 54.66\MBunit &
 0.01\secunit & 1.95\MBunit \\
 mult6x6 & 117 & 106 &
 0.03\secunit & 54.92\MBunit &
 0.03\secunit & 1.95\MBunit \\
 mult8x8 & 203 & 188 &
 0.40\secunit & 55.43\MBunit &
 0.96\secunit & 2.21\MBunit \\
 mult10x10 & 313 & 294 &
 18.11\secunit & 85.91\MBunit &
 22.85\secunit & 3.61\MBunit
\end{tabular}
}
\end{center}
\caption{Deciding satisfiability with \PolyBoRi using Gr\"{o}bner basis computations in comparison with MiniSat, a state-of-the-art SAT solver.}
\label{tab:pb_sat}
\end{table}
Although the memory consumption of \PolyBoRi is larger, \tabref{tab:pb_sat}
illustrates that
the computation time of both approaches is comparable for this kind of practical
examples.
(The first part of the table  was computed using the preprocessing
motivated by \thmref{unique-generator}.)
In  particular, it shows, that in our research area the algebraic
approach is competitive with SAT-solvers.

The advantages of \PolyBoRi are illustrated by the examples above as follows:
the fast Boolean multiplication can be seen in the
pigeon hole benchmarks.
%
The computations of the uuf problems include a large number
of generators, consisting of initially short polynomials, which
lead to large intermediate results.
The algorithmic improvement of \symmgbGFTwo and the \optimised pair handling
render the treatment of these example with algebraic methods possible.


In this way the initial performance of \PolyBoRi is promising. The
data show that the advantage of \PolyBoRi grows with the number of
variables. For many practical applications this size will be even bigger.
Hence, there is a chance,
that it will be possible to tackle some of these problems in
future by using more \specialised approaches. A key point in the
development of \PolyBoRi is to facilitate problem specific and high performance
solutions.

\section{Conclusions}

For efficient treatment of bit-level formulations of digital systems
we have developed
\specialised methods for the analysis of polynomial systems in Boolean rings,
that is quotient rings of the form~$\Ztwoxoneton$ modulo the field polynomials.
For this purpose improvements were achieved on multiple levels.
On one hand, a tailored data structure was introduced to represent Boolean
polynomials which correspond to
canonical representatives of the elements in the quotient ring.
This structure, which is derived from
zero-suppressed binary decision diagrams~(ZDDs),
is compact and allows to
apply operations used in \Groebner basis computations in reasonable time.
Further,
enhancement were due to the \specialised 
\Groebner basis  algorithm
\symmgbGFTwo itself.
Exploiting special properties in the Boolean case, special
criteria for keeping the set of critical pairs small were proposed.
In addition,  (recursive) caching of previous computations and \utilising
symmetry makes it possible to efficiently reuse  results arising from likewise
polynomials.
Also, the \PolyBoRi system, a framework for Boolean rings, was presented as
reference implementation for \symmgbGFTwo and for the ZDD-based data
structure representing Boolean polynomials.

Word-level formulations of digital systems lead us to
investigate \Groebner bases over rings. More generally, we developed the
theory of standard bases over rings for which systems of linear equations can be
solved effectively. For weak factorial principal ideal rings
we developed special criteria for $s$-polynomials and for the normal form
algorithm which proved effective.

%
The \PolyBoRi framework for Boolean \Groebner bases showed that~-- in particular
if  there are no immediate counter examples~-- the proposed approach
has already reached the same level as a state-of-the-art SAT-solver at least
for some standard benchmark examples. The
advantage of an effective theory of Boolean Gr\"{o}bner basis is, that they are
a general and flexible tool which opens the door to computational algebra over Boolean rings.

\bibliography{groebner,pbori,groebner-oli}
\end{document}